\documentclass[11pt]{amsart}

\usepackage{a4wide}
\usepackage[T1]{fontenc}
\usepackage[latin1]{inputenc}
\usepackage{amssymb}
\usepackage{amsfonts}
\usepackage{latexsym}
\usepackage{psfrag}
\usepackage{epsfig}
\makeatletter
\newcommand{\cal}{\mathcal}
\newcommand{\R}{\mathbb{R}}
\newcommand{\C}{\mathbb{C}}
\newcommand{\bbS}{\mathbb{S}}
\newcommand{\bbT}{\mathbb{T}}

\newcommand{\calS}{\cal {S}}

\newcommand{\N}{\mathbb{N}}

\newcommand{\al}{\alpha}
\newcommand{\be}{\beta}
\newcommand{\si}{\sigma}

\newcommand{\Om}{\Omega}
\newcommand{\om}{\omega}

\newcommand{\ti}{\tilde}
\newcommand{\Ti}{\widetilde}

\newcommand{\To}{\longrightarrow}

\newcommand{\de}{\delta}

\newcommand{\im}{\mbox{Im}}
\newcommand{\re}{\mbox{Re}}

\newcommand{\omout}{\om_{\rm{out}}}
\newcommand{\omin}{\om_{\rm{in}}}
\newcommand{\sgn}{\mbox{sgn}}

\def \.{{\bf{\cdot}}}

\def \dx{{\frac{d}{dx}}}
\def \dz{{\frac{d}{dz}}}
\def \ev{{\mbox{\begin{footnotesize}even\end{footnotesize}}}}
\def \od{{\mbox{\begin{footnotesize}odd\end{footnotesize}}}}
\def \im{{\mbox{Im }}}
\def \re{{\mbox{Re }}}
\def \W{{\cal W}}
\newtheorem{prop}{Proposition}[section]
\newtheorem{lem}{Lemma}[section]
\newtheorem{thm}{Theorem}[section]
\newtheorem{rem}{Remark}[section]

\newtheorem{defi}{Definition}[section]

\numberwithin{equation}{section}

\begin{document}

\makeatother
\title[Klein paradox and Scattering theory]{Klein paradox and Scattering theory for the semi-classical 
 Dirac equation 
\vspace{5mm}}
\author{Abdallah Khochman}
\date{\today}
\email{Abdallah.Khochman@math.u-bordeaux1.fr}
\address{Universit\'e Bordeaux I, Institut de Math\'ematiques, UMR CNRS 5251, 351, cours de la Lib\'eration, 33405 Talence, France}
\maketitle{}


\begin{abstract}
We study the Klein paradox for the semi-classical 
 Dirac operator on $\R$ with potentials having constant limits, not necessarily the same at infinity. Using the
complex WKB method, the time-independent scattering
theory in terms of incoming and outgoing plane wave solutions is
established. The corresponding scattering matrix is unitary. We obtain an 
 asymptotic expansion, with respect to the semi-classical  parameter $h$,
 of the scattering matrix in the cases of the Klein paradox, the total transmission and the total reflection.
  Finally, we treat the scattering problem in the zero mass case.
\end{abstract}

\maketitle{{\it Keywords}: Semi-classical Dirac operator - Scattering matrix - Klein paradox - Complex WKB method.

{\it Mathematics classification:} 81Q05 - 47A40 - 34L40 - 34E20 - 34M60. }

\section{Introduction}
In mathematics and physics, the scattering theory is a framework for studying and understanding the scattering of waves and particles.
The scattering matrix for the one-dimensional Dirac operator $H$ is closely related to the transition probability of particles through
a potential. However, if the potential does not vanish at infinity, a Klein paradox might occur. The latter is of great historical importance
 in order to justify the existence of the antiparticle of an electron (the positron) and explaining qualitatively the pair creation process
 in the collision of particle beam with strongly repulsive electric field. The explanation of this Klein paradox usually resorted to the
 concept of "hole" in the "negative-energy electron sea".
For more physical interpretations we refer to  Klein \cite{OK}, Sauter \cite{FS}, Bjorken-Drell \cite{JBSD}, Sakurai \cite{JS},
Thaller \cite{BT} and Calogeracos-Dombey \cite{ACND} for the history of the Klein paradox. This paradox appears also for the
 Klein-Gordon equation, here no concept of "hole" is needed (see Winter \cite{RGW} and Ni-Zhou-Yan \cite{GNWZJY} for a constant
  potential at infinity and Bachelot \cite{AB} for an electrostatic potential having different asymptotics at infinity).
  The comparison between the Klein paradox for this two equations 
has been discussed in \cite[Part C]{RGW}. A Klein paradox phenomenon occurs also in quantum field theory (see Hund \cite{FH}
 and Manogue \cite{CAM}). It is clear that this paradox cannot appear  for Schr\"odinger operators.

\medskip
For a scalar potential having real limits $V^\pm$ at $\pm\infty$,
the Klein paradox of the Dirac equation occurs if $V^+-V^->2mc^2$.
In this case the higher part of $\si(H)$ intersects its lower part.
If the energy $E$ is in this intersection, for a wave-packet which
comes from the left and moves towards the potential,
 a part of it is reflected, another part being transmitted.
 The transmitted part moves to the right and behaves like a solution with negative energy. 
 Ruijsenaars-Bongaarts \cite{SRPB} (see also Thaller \cite{BT}) have mathematically treated the Klein paradox and the
 scattering theory for the Dirac equation with one-dimensional potentials constant outside a compact set.
 They have established the connection between time-dependent and
time-independent scattering theory in terms of incoming and outgoing plane wave
solutions. The exact calculus of the scattering matrix for one-dimensional Dirac operator is only
known for a few number of explicit potentials (see Klein \cite{OK} for a rectangular step potential and Fl\"ugge \cite{SF} for the potential $V(x)=\mbox{tanh}(x)$).
Nevertheless, we are neither aware of works dealing with the asymptotic expansion of the scattering matrix,
with respect to the semi-classical parameter $h$.

\medskip
For one-dimensional Schr\"odinger operators, 
there are several approaches which have been developed dealing with
the computation of the transmission coefficient through a barrier.
Ecalle \cite{JE} and Voros \cite{AV} have developed the so-called
complex WKB analysis which gives approximations in the complex plane
of the solutions of a Schr\"odinger equation. This approach is used
in a new formalism by Grigis for the Hill's equation \cite{AG}. This
method is also used by C. G\'erard-Grigis \cite{CGAG} to calculate
the eigenvalues near a potential barrier and by Ramond \cite{TR} for
scattering problems.
 For references and a historical discussion, we refer to Ramond \cite{TR}. The complex WKB
 method has been extended to 
 a class of Schr\"odinger systems by
Fujii\'e-Lasser-N\'edelec \cite{SFCLLN}.

\medskip
%
%
The purpose of this paper is to give an asymptotic expansion, with
respect to the semi-classical parameter $h$, of the coefficients of
the scattering matrix for the one-dimensional Dirac operator with
potentials having different limits at infinity. We establish the
exponential decay of the transmission coefficient in the Klein
paradox case (cf. Theorem \ref{thm2point} below).
  We calculate the coefficients of the scattering matrix in terms of incoming and outgoing
  solutions. Therefore, we use the complex
WKB analysis to construct solutions of the Dirac equation. The
usefulness of
this analysis is that it provides, rather than approximate solutions with error bounds, solutions in the complex plane with a complete asymptotic expansion with respect to the semi-classical parameter $h$
, with a priory estimates on the coefficients.

\medskip
The paper is organized as follows. In the next section, 
 we introduce the perturbed Dirac operator on $\R$, study the time-independent scattering theory 
 and state our main 
 results. In Section \ref{sectionconst}, we develop the complex WKB method 
 and show a
complete asymptotic expansion 
 of the
coefficients. 
 In Section \ref{sectionJost}, the existence of incoming and
outgoing Jost solutions is proved.
   In
Section \ref{sect-2ptourn}, we analyze the semi-classical behavior of
the scattering matrix in the Klein paradox case. 
 The total transmission over a potential
barrier and the total reflection 
 are
studied in Section \ref{sec-0ptourn} and Section \ref{sec-1ptourn}. Finally, in
Section \ref{sec-mass0}, we study the Klein paradox for the zero mass case.
\section{Assumptions and results}\label{sectionassuptionresutls}
We consider the self-adjoint Hamiltonian $H=H_0+V$, where $H_0$ is
the semi-classical Dirac operator on $\R$:
\begin{eqnarray}\label{opDirac}
H_0=-ihc\al\frac{d}{dx}+mc^2\be,
\end{eqnarray}
with domain $D(H_0)=H^1(\R)\otimes\C^2\subset{\cal
H}=L^2(\R)\otimes\C^2$, where $h\searrow0$ is the semi-classical
parameter, $m\geq0$ is the mass of the Dirac particle and$\;c$ is
the celerity of the light. The coefficients $\al,\be$ are the
$2\times2$ Pauli matrices satisfying the anti-commutation relation
$$\al\be+\be\al=0,$$ and $\al^2=\be^2=I_2,$ where $I_2$ is the $2\times2$
identity matrix.

The operator $V$ is the multiplication by $VI_2,$ where $V$ is a smooth 
 electrostatic potential satisfying:
\\
\\
(A):   {\it The function $V$ is real on the real axis, analytic in the sector}
$$\calS=\{x\in\C,\ |\mbox{Im }x|<\epsilon|\mbox{Re }x|+\eta\}
,$$
{\it for some $\epsilon,\eta
>0$, and satisfies the following estimates:}
\begin{eqnarray}\label{eqpotential}
|
(V(x)-V^{\pm}
)|=O(\langle x\rangle^{-\de})\ \ \ for\  \re (x)\To\pm\infty\ \  in\
\ \calS.
\end{eqnarray}
$Here$, $\langle x\rangle=(1+|x|^2)^\frac12$, $\de>1\ $ $and$ $\ V^-<V^+$.
{\begin{center}
\psfrag{0}[rh][][1][0]{$0$}
\psfrag{Vl}[][][1][0]{$V^-$}
\psfrag{Vr}[][][1][0]{$V^+$}
\psfrag{titre}[][][1][0]{Fig. 1. The potential $V$}
\includegraphics[height=5cm,width=15cm]{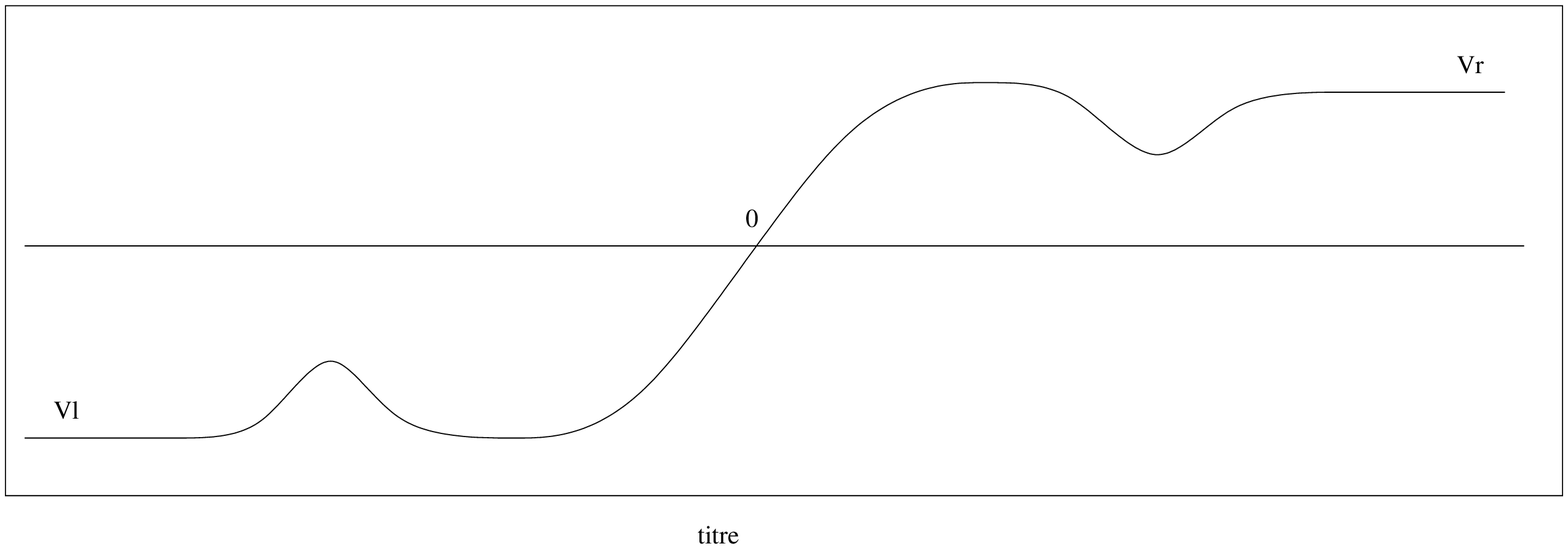}
\end{center}}

The spectrum of the free Dirac operator $H_0$ is
$
]-\infty,-mc^2]\cup[mc^2,+\infty[$ and it 
 is purely absolutely continuous. Under assumption  (A) 
 the operator $H=H_0+V$ is a self-adjoint operator 
 and has essential spectrum  (see Appendix \ref{spectrumofdiracoperator}):\begin{eqnarray}\label{sepectreess}
\si_{ess}(H)&=&
]-\infty,-mc^2+V^+]\cup[mc^2+V^-,+\infty[.\end{eqnarray}

There are several representations of the matrices $\al,\be.$ For
example, Hiller \cite{JH} used $\al=\si_2,\;\be=\si_3,$ Nogami and
Toyoma \cite{YNFT} used $\al=\si_2,\;\be=\si_1$, where
$\si_j$, $j=1,2,3,$ are the standard representation for Dirac-Pauli
matrices. Most calculations with Dirac matrices can be done without referring
to a particular representation 
$($see Thaller \cite[Appendix 1A]{BT}$)$. 
Here, we choose the $1+1$ dimensional representation of the Dirac matrices
\begin{eqnarray}\label{eqrep}
\al=\si_1=\left(\begin{array}{cc}
0&1\\
1&0
\end{array}\right),\;\;\be=\si_3=\left(\begin{array}{cc}
1&0\\
0&-1
\end{array}\right).
\end{eqnarray}

The solutions of
\begin{eqnarray}\label{eq valeu propre}
Hu=\left(-ihc\si_1\frac{d}{dx}+mc^2\si_3+V(x)I_2\right)u=Eu,\ \ \ \ E\in\R,
\end{eqnarray}
should behave as $x\To \pm\infty$ like
$$a_{\pm}^+(E,h)\exp(+\frac{1}{hc}(m^2c^4-(V^{\pm}-E)^2)^{1/2}x)+a_{\pm}^-(E,h)\exp(-\frac{1}{hc}(m^2c^4-(V^{\pm}-E)^2)^{1/2}x).$$
Here, the square root $(\cdot)^\frac12$ is to be defined more precisely according to the sign of $m^2c^4-(V^{\pm}-E)^2$.

In the following, we will use these 
 intervals on the E-axis:

$\ \ $ I.  $\ $ $V^++mc^2<E$,

$\ $ II. $\ $ max($mc^2+V^-,V^+-mc^2$)$<E<V^++mc^2$,

$\ $III. $\ $ $V^-+mc^2<E<V^+-mc^2\ \ \ $ if $\ \ \ V^+-V^->2mc^2$,

$\ $IV. $\ $ $V^--mc^2<E<\ $min($V^-+mc^2,V^+-mc^2$),

$\  $ V.  $\ $ $E<V^--mc^2$.\\
If $ V^+-V^->2mc^2$, the different regions are represented in the following figure: 
{\begin{center}
\psfrag{0}[rh][][1][0]{$0$}
\psfrag{Vr+mc2}[][][1][0]{$V^++mc^2$}
\psfrag{Vr-mc2}[][][1][0]{$V^+-mc^2$}
\psfrag{Vl+mc2}[][][1][0]{$V^-+mc^2$}
\psfrag{Vl-mc2}[][][1][0]{$V^--mc^2$}
\psfrag{V(x)}[][][1][0]{$\ \ E
$} \psfrag{x}[][][1][0]{$ $} \psfrag{regionI}[][][1][0]{Region I}
\psfrag{regionII}[][][1][0]{Region II}
\psfrag{regionIII}[][][1][0]{Region III}
\psfrag{regionIV}[][][1][0]{Region IV}
\psfrag{regionV}[][][1][0]{Region V} \psfrag{titre}[][][1][0]{Fig.
2. Different regions on the E-axis}
\includegraphics[height=6cm,width=15cm]{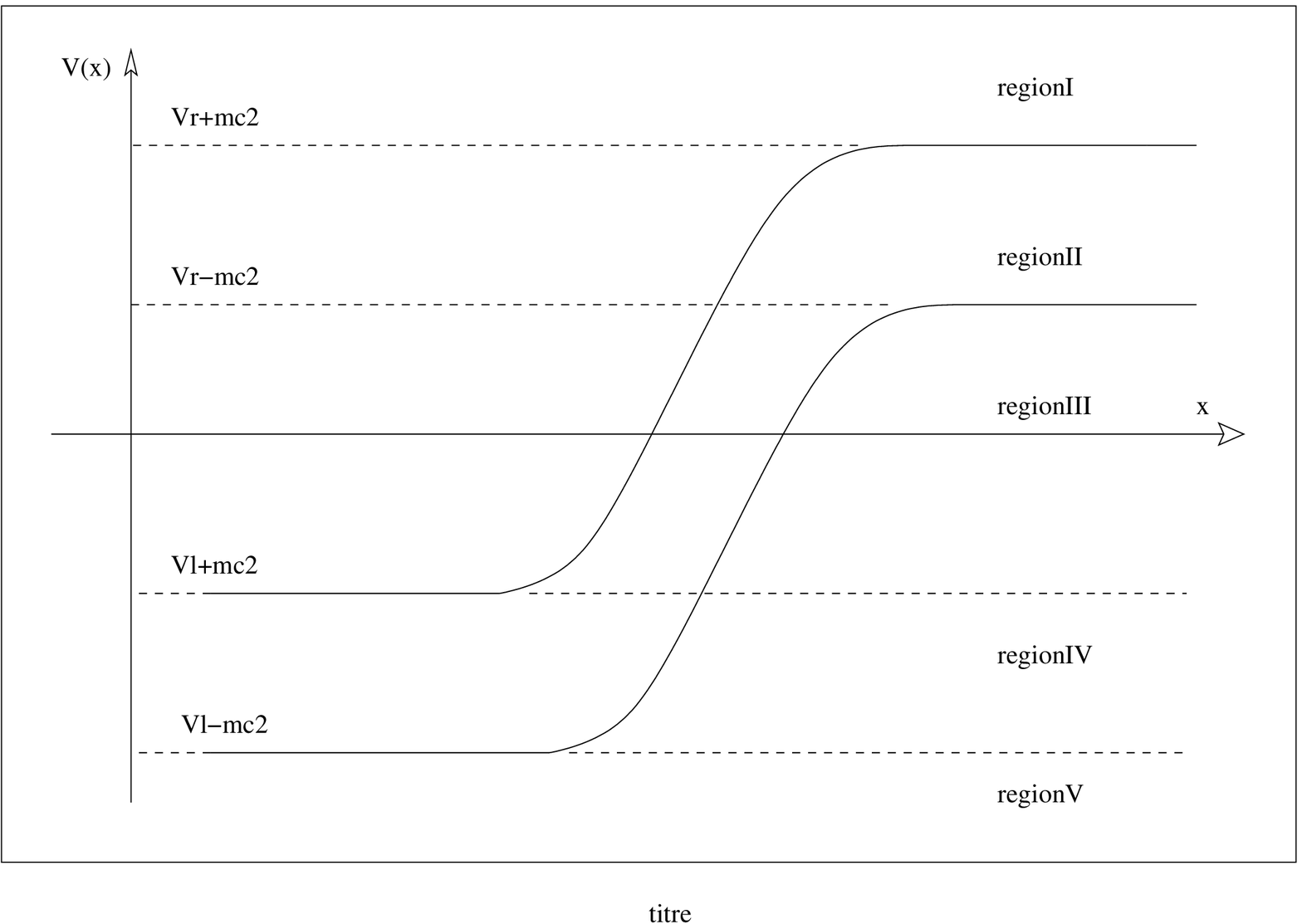}
\end{center}}
We study the semi-classical behavior of the scattering matrix for 
 the different values of the
energy $E$. Let us describe now the time-independent scattering
problem briefly.
 For $E\in$ I, III or  V, the four Jost solutions $\omin^\pm$, $\omout^\pm$
 (see Theorem \ref{solutioJost}) are the solutions of (\ref{eq valeu propre}) which behave exactly as
\begin{eqnarray}
\label{ominpmintro} \omin^\pm&\sim&\exp\{\mp\frac{i}{hc}\Phi(E-V^\pm)
x\}\left(\begin{array}{c}
A(E-V^\pm)\\
{\mp}{A(E-V^\pm)^{-1}}
\end{array}\right)\,\;\;\;\mbox{as }x\To\pm\infty,\\
\label{omoutintro} \omout^\pm&\sim&\exp\{\pm\frac{i}{hc}\Phi(E-V^\pm)
x\}\left(\begin{array}{c}
A(E-V^\pm)\\
\pm{A(E-V^\pm)^{-1}}
\end{array}\right)\,\;\;\;\mbox{as }x\To\pm\infty,
\end{eqnarray}
with $\Phi(E)=\sgn(E)\sqrt{E^2-m^2c^4}$, $A(E)=\sqrt[4]{\frac{E+mc^2}{E-mc^2}}$  and $\sgn(E)=\frac{E}{|E|}$ for $E\not\in[-mc^2,\,mc^2]$. 
 Analogous definitions of Jost solutions can be found in the works of Ruijsenaars-Bongaarts \cite{SRPB} and Thaller \cite{BT} for one-dimensional step potentials. In this paper, we denote $\sqrt{x},\ \sqrt[4]{x}$ the positive determination of $x\in\R^+\To (x)^{\frac{1}{2}},\ (x)^{\frac{1}{4}}$ respectively.

The 
ordinary scattering problem is the following: what are the
components  of a solution $u$ of the Dirac equation (\ref{eq valeu propre}) in the basis
($\omout^+,\ \omout^-$) of the outgoing Jost solutions, knowing its
component in the basis ($\omin^-,\ \omin^+$) of the incoming Jost
solutions. The $2\times 2$ matrix relating these coefficients is
called the scattering matrix and we will denote it by $$
{\bbS}=\left(\begin{array}{cc}
 s_{11}& s_{12} \\
s_{21} & s_{22}
 \end{array}\right)
.$$
Precisely, if we take $u$ 
 a solution of (\ref{eq valeu propre}), $$u=a_{\rm{in}}\omin^-+b_{\rm{in}}\omin^+= a_{\rm{out}}\omout^++b_{\rm{out}} \omout^-,$$
 the scattering matrix is such that
$${\bbS}\left(
\begin{array}{l}
 a_{\rm{in}}\\
b_{\rm{in}}
        \end{array}\right)=\left(\begin{array}{l}
a_{\rm{out}} \\
b_{\rm{out}}
        \end{array}\right),
$$
which is equivalent to
\begin{eqnarray}\label{matriceSoutin}
(\omin^-,\omin^+)=(\omout^+, \omout^-){\bbS}.
\end{eqnarray}
 Since $V$ is real on  the real axis,  we have (see (\ref{eqrep}))
  \begin{eqnarray}\label{relation-entre soljost}
 \overline{\omin^\pm}&=&\be\omout^\pm.
 \end{eqnarray}
We also have the following relations between the coefficients of $\bbS(E,h)$:
\begin{eqnarray}\label{relationmatricescattering}
s_{11}(E,h)=s_{22}(E,h)\ \ \ \ \mbox{and}\ \ \ \ s_{12}(E,h)=-{\overline {s_{21}}}(E,h)\frac{s_{11}(E,h)}{{\overline {s_{11}}}(E,h)},
\end{eqnarray}
so that $s_{11}$ and $s_{12}$ determine completely the scattering matrix.\\


The reflection and transmission coefficients $R(E,h)$ and
$T(E,h)$ are, by definition, the square of the  modulus of the
coefficients $s_{21}$ and $s_{11}$ respectively. They correspond to
the probability for a purely incoming-from-the-left particle to be
reflected to the left or transmitted to the right. Using (\ref{relation-entre soljost}), (\ref{relationmatricescattering}) and calculating the determinant of (\ref{matriceSoutin}), we have the well-known relation $R(E,h)+T(E,h)=1$ and, the scattering matrix $\bbS(E,h)$ is unitary.

To calculate the scattering matrix $\bbS(E,h)$ we will use the transfer matrix $\bbT$, which is defined by
$$(\omin^-,\ \omout^- )=(\omout^+,\ \omin^+)\bbT.$$

The determinant of this matrix is equal to $1$ since the two Wronskians $\W( \omin^-,\ \omout^-)$ and  $\W( \omout^+,\ \omin^+)$ are equal to $-2$ (see Definition \ref{defwronskian}). Using the relation
 (\ref{relation-entre soljost}), we obtain that $\bbT$ is determined by two coefficients:
\begin{eqnarray}\label{matrixT}{\bbT}=\left(\begin{array}{cc}
 t(E,h)& r(E,h)\\
 \overline{\,r\,}(E,h)&\overline{\,t\,}(E,h)
 \end{array}\right).
\end{eqnarray}
Moreover, using that $\det({\bbT})=1$, we obtain
\begin{eqnarray}\label{det=1}
|t(E,h)|^2-|r(E,h)|^2=1.
\end{eqnarray}
Consequently, we can write the scattering matrix in terms of the coefficients of the transfer matrix $\bbT$:
\begin{eqnarray}\label{matrixS}{\bbS}=\frac{1}{\overline{\,t\,}(E,h)}\left(\begin{array}{cc}
 1&- r(E,h)\\
 \overline{\,r\,}(E,h)&1
 \end{array}\right).
\end{eqnarray}

We will use WKB approaches to describe the amplitude of the
coefficients of the scattering matrix for $h\searrow0$. For these,
let us introduce the following definition.
\begin{defi}\label{CAS}(See Sj\"ostrand \cite{JSjostrand})
A function $f(z,h)$ defined in $U\times]0,h_0[$, where $U$ is an
open set in $\C$ and $h_0>0$, is called a classical analytic symbol
(CAS) of order $m\in\N$ in $h$ if $f$ is an analytic function of
$z\in U$ and if there exists a sequence $(a_j(z))$ of analytic
functions in $U$ such that
\begin{itemize}
\item For all compact set $K\subset U$, there exists $C>0$ such that, for all $z\in K$, one has $$|a_j(z)|\leq C^{j+1}j^j.$$
\item 
The function $f(z,h)$ admits the series $\sum_{0\leq j\leq +\infty 
}a_j(z)h^{m+j}$ as asymptotic expansion as $h$ goes to zero in the following sense. 
For any $C_1>C$, 
 we have $$f(z,h)-\sum_{0\leq j\leq h^{-1}/eC_1}a_j(z)h^{m+j}=
O(e^{-\rho/h}),$$ for some $\rho>0$ and all $z\in K$\end{itemize}
\end{defi}
The main theorem concerning the Klein paradox case for $m>0$ (i.e. for the energy level $E\in$ III) is the following: 
\\
{\begin{center}
\psfrag{E}[][][1][0]{E}
\psfrag{t1}[rh][][1][0]{$\!\!t_1(E)$}
\psfrag{t2}[][][1][0]{$\ t_2(E)$}
\psfrag{Vl+mc2}[][][1][0]{$V^-+mc^2$}
\psfrag{Vl-mc2}[][][1][0]{$V^--mc^2$}
\psfrag{Vr+mc2}[][][1][0]{$V^++mc^2$}
\psfrag{Vr-mc2}[][][1][0]{$V^+-mc^2$} \psfrag{titre}[][][1][0]{Fig.
3. Graph of $V(x)+mc^2$ and $V(x)-mc^2$}
 \hspace{2cm}\includegraphics[height=6cm, width=17cm]{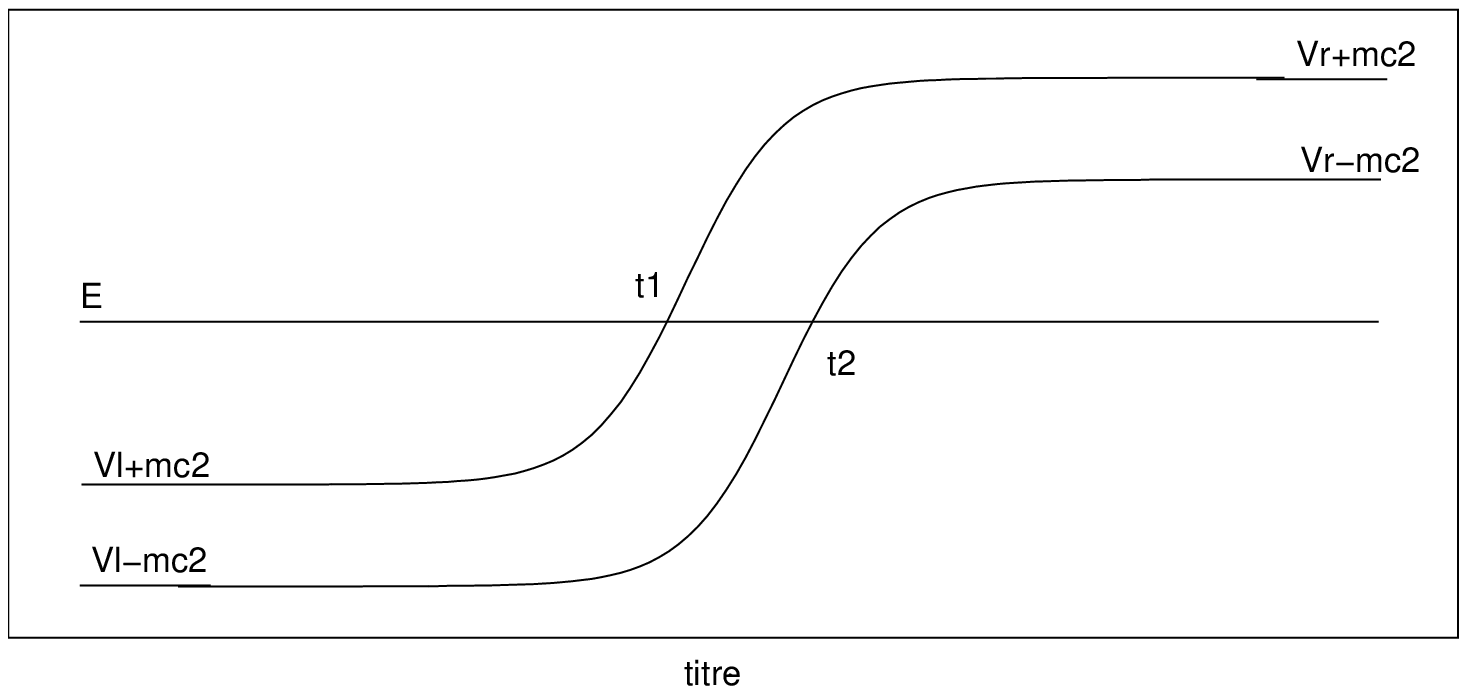}
\end{center}
}
\begin{thm}\label{thm2point}\emph{[}{\bf Klein paradox}\emph{]}
Let $V$ be a potential satisfying assumption (A), $E\in$ \emph{III}
and $m>0$. Suppose that there are only two  simple zeros
$t_1(E)<t_2(E)$ of $m^2c^4-(V(x)-E)^2$ (see Fig. 3). Then there
exists three classical analytic symbols $\phi_1(h)$, $\phi_2(h)$ and
$\phi_3(h)$ of non-negative order such that:
\begin{eqnarray}
s_{11}&=&s_{22}=(1+h\phi_1(h))\exp\{-S(E)/h\}\exp\{iT(E)/h\},
\label{eqs11mnonul}\\\nonumber\\
s_{21}&=&(i+h\phi_2(h))\exp\left\{\frac{2i}{h}\left(t_1(E)\sqrt{E^-}+\int_{-\infty}^{t_1(E)}Q^-(t,E)dt\right) \right\},\label{eqs21mnonul} \\\nonumber\\
s_{12}&=&(i+h\phi_3(h))\exp\left\{\frac{2i}{h}\left(t_2(E)\sqrt{E^+}-\int_{t_2(E)}^{+\infty}Q^+(t,E)dt\right) \right\},\label{eqs12mnonul}
\end{eqnarray}
where $S(E)$ is the classical action between the two turning points $t_1(E)$ and $t_2(E)$
$$S(E)=\int_{t_1(E)}^{t_2(E)}\sqrt{\frac{m^2c^4-(V(t)-E)^2}{c^2}}dt.$$
Moreover
$$ Q^{-}(t,E)=\sqrt{\frac{-m^2c^4+(V(t)-E)^2}{c^2}}-\sqrt{E^{-}}, \ \ \ \emph{for }t<t_1(E),$$
$$ Q^{+}(t,E)=\sqrt{\frac{-m^2c^4+(V(t)-E)^2}{c^2}}-\sqrt{E^{+}}, \ \ \ \emph{for }t>t_2(E),$$
and
$$T(E)=\int_{-\infty}^{t_1(E)}Q^-(t,E)dt-\int_{t_2(E)}^{+\infty}Q^+(t,E)dt+t_1(E)\sqrt{E^-}+t_2(E)\sqrt{E^+},$$
where $$E^{\pm}=\frac{-m^2c^4+(V^{\pm}-E)^2}{c^2}.$$
\end{thm}
We remark that this scattering matrix behaves like in the case of
the Schr\"odinger operator with a barrier potential. In particular
the term
$e^{-S(E)/h}$ which decays exponentially, can be viewed as a tunneling
effect (see Ramond \cite[Theorem 1]{TR}).\\\\

In the zero mass case, we have the following theorem:
\begin{thm}\emph{[{\bf Zero mass case}]}\label{thmmassnull} 
Let $V$ be a potential satisfying assumption (A), $E\in$ \emph{III} and $m=0$. Suppose that there is only a simple zero $t_0(E)$ of $V(x)-E$. Then,
there is a classical analytic symbol $\phi(h)$ 
such that:
\begin{eqnarray}
s_{11}&=&s_{22}=(1+h\phi(h))\exp\{iT_0(E)/h\},
\label{eqs11mnonul}\\\nonumber\\
s_{21}&=&O(h)
,\\\nonumber\\
s_{12}&=&O(h)
\label{eqs12mnul}.
\end{eqnarray}
 Here,
\begin{eqnarray*}
T_0(E)&=&T(E)\Big|_{m=0}\\ &=&\frac{1}{c}\left(\int_{-\infty}^{t_0(E)}(-V(t)+V^-)dt-\int_{t_0(E)}^{+\infty}(V(t)-V^+)dt+t_0(E)(V^+-V^-)\right).
\end{eqnarray*}
\end{thm}
\begin{rem}
We can not permute the limits of the scattering matrix $\bbS$ as $m\to0$ and $h\to0$.
Indeed, if we take the limits of $s_{12}$ in (\ref{eqs12mnonul}) and
(\ref{eqs12mnul}), we obtain
$$ \lim_{m\to0}\lim_{h\to0}|s_{12}|=1,\ \ \ \ \lim_{h\to0}\lim_{m\to0}|s_{12}|=0.$$
\end{rem}
Now, we come back to 
 the non-zero mass case and we treat reflection and transmission cases (see Sections \ref{sec-1ptourn},
\ref{sec-0ptourn}).\\

If we take the energy level $E\in$II, there are two Jost solutions $\omin^-,\ \omout^-$ satisfying (\ref{ominpmintro})
and (\ref{omoutintro}) for $x\To-\infty$ and 
there does not exist an oscillating solution for $x\To+\infty$.
Instead, as $x\to +\infty,$ there exists an exponentially decaying
solution and an exponentially growing solution. Since the last
function doesn't represent a physical state we limit ourself to the
one dimensional space generated by the decaying solution $\om_d^+$
(unique up to a constant). This function satisfies (see Theorem
\ref{solutioJost}):
\begin{eqnarray}
 \om_d^+&\sim&\exp\{-\frac{1}{hc}\sqrt{m^2c^4-(V^+-E)^2}x\}\left(\begin{array}{c}
- i\sqrt[4]{\frac{mc^2+E-V^+}{mc^2-E+V^+}}\\
\sqrt[4]{\frac{mc^2-E+V^+}{mc^2+E-V^+}}
\end{array}\right)\,\;\;\;\mbox{as }x\To+\infty.
\end{eqnarray}
In this case we have
\begin{thm}\label{thmt1pt}\emph{[{\bf Total reflection}]}
Let $V$ be a potential satisfying assumption (A), $E\in$ \emph{II}
 and $m>0$. Suppose that there is only a simple zero $t_1(E)$ of $m^2c^4-(V(x)-E)^2$. Then the vector space of the solutions of
$(H-E)u=0$ with $u$ bounded 
 is a one
dimensional space generated by
$$u=\omin^-+\al_{\rm out}^- \omout^-,\ \ \ \mbox{with}$$
\begin{eqnarray}\label{valeuralphal-}
\al_{\rm
out}^-=-i(1+h\phi_1(h))\exp\left\{\frac{2i}{h}\left(\int_{-\infty}^{t_1(E)}Q^-(t,E)dt
+\sqrt{E^-}t_1(E)\right)\right\}.
\end{eqnarray}
Moreover
$$u=\be_d^+ \om_d^+,\ \ \ \ \mbox{with}$$
 \begin{eqnarray*}\be_d^+&=&e^{i\pi/4}(1+h\phi_2(h))\times\\& &\exp\left\{\frac1h\left(\int_{+\infty}^{t_1(E)}Q_-^+(t,E)dt
 +\sqrt{-E^+}t_1(E)+i\int_{-\infty}^{t_1(E)}Q^-(t,E)dt+i\sqrt{E^-}t_1(E)\right)\right\}.\end{eqnarray*}
Here
$$Q_-^+(t,E)=\sqrt{\frac{m^2c^4-(V(t)-E)^2}{c^2}}-\sqrt{-E^+},$$
$Q^-(t,E)$ and $E^\pm$ are the functions of Theorem \ref{thm2point} 
 and $\phi_j(h)$, $j=1,2$ are classical analytic symbols of non-negative order.
\end{thm}
For $E\in$ IV, there is also total reflection cases which can treated similarly to the previous theorem. As in \cite{SRPB},
there is also a scattering interpretation of the
previous theorem. Since we work in a one-dimensional space, the
scattering matrix is now a scalar.
\begin{rem}\label{Scattering interpretation}\emph{[{\bf Scattering interpretation}]}
 We call $u^{\rm{in}}=\omin^-+\al_{\rm
out}^-\omout^-$ the ``$incoming$'' solution. In the same way, there
exists  a unique bounded solution
$$u^{\rm{out}}=\omout^-+\al_{\rm in}^-\omin^-,$$
 which is called the ``$outgoing$'' solution.

If $u$ is a bounded solution of $(H-E)u=0$ (i.e. $u=Au^{\rm{in}}$)
then $u=Bu^{\rm{out}}$. The scattering matrix $\bbS$ is defined by 
$$B=\bbS A.$$
From (\ref{valeuralphal-}), we have
$$\bbS=\al_{\rm out}^-=-i(1+h\phi(h))\exp\left\{\frac{2i}{h}\left(\int_{-\infty}^{t_1(E)}Q^-(t,E)dt
+\sqrt{E^-}t_1(E)\right)\right\}
.$$
\end{rem}

 For $E\in$ I or V, a total transmission phenomena occur:
\begin{thm}\emph{[{\bf Total transmission}]}\label{thmtransmission}
Let $V$ be a potential satisfying assumption (A), $E\in$ \emph{I},
$m\geq0$ and $m^2c^4-(V(x)-E)^2\neq0$. 
Then there are a classical analytic symbol $\phi(h)$ and positive
constant $C$
 such that:
\begin{eqnarray}
s_{11}&=&s_{22}=(1+h\phi(h))\exp\{i\Ti T(E)/h\},
\\\nonumber\\
s_{21}&=&O(e^{-C/h})
\ \ \ \ \mbox{and    }\ \ \ s_{12}=O(e^{-C/h}),
\end{eqnarray}
where 
$$\Ti T(E)=\int_{-\infty}^{0}Q^-(t,E)dt+\int_{0}^{+\infty}Q^+(t,E)dt,$$
and $Q^-(t,E),\ Q^+(t,E)$ are the functions of Theorem \ref{thm2point} defined here for any $t\in\R$.
\end{thm}
We can calculate the scattering matrix for $E\in$ V in the same way
of $E\in$ I. We remark that the behavior of the incoming and
outgoing Jost solutions exchanges between these two cases. This is
in agreement with the physical interpretation (see
\cite[p.121]{BT}).


\section{Complex WKB solutions}\label{sectionconst}
We wish to find a representation formula for the solutions of (\ref{eq valeu propre}),  
 from which it is possible to deduce the asymptotic expansion in $h$. The method is known as complex WKB method.
See \cite{TR} \cite{CGAG}, \cite{SFTR1}, \cite{SFTR2}, \cite{SFCLLN} for constructions of solutions of the Schr\"odinger equation.\\ \\
In a complex domain $
\calS$, we study the Dirac system (\ref{eq valeu propre}
)
  which is of the form
\begin{eqnarray}\label{eqDirac1}
(H-E)u(x)=\left(\begin{array}{cc}
mc^2+V(x)-E&-ihc\dx\\
-ihc\dx&-mc^2+V(x)-E
\end{array}\right)u(x)=0,
\end{eqnarray}
or equivalently
\begin{eqnarray}\label{eqDirac2}
\frac hi\dx v(x)=\left(\begin{array}{cc}
0&g_+(x)\\
-g_-(x)&0
\end{array}\right)v(x),
\end{eqnarray}
where $v(x)=\left(\begin{array}{cc}
0&1\\
1&0
\end{array}\right)u(x)=M^{-1}u,$ and the functions $$g_{\pm}(x)=\frac{-mc^2\mp(V(x)-E)}{c},$$ are holomorphic in $\calS
$. The following considerations will lead to the construction of complex WKB solutions for Dirac system.
\subsection{Formal construction}
First, we introduce a new complex coordinate
\begin{eqnarray}\label{eqchange}
z(x)=z(x,x_0)=\int_{\gamma(x_0,x)}(g_+(t)g_-(t))^\frac12dt=\int_{x_0}^x(g_+(t)g_-(t))^\frac12dt,\;\;\;\;x_0\in D.
\end{eqnarray}
One of our tasks will be of course to choose 
 the simply connected subset  $D$ of
$\calS$ 
such that $t\To(g_+(t)g_-(t))^\frac12$ 
 is well-defined,  but let's work
formally for a while. The $\gamma(x_0,x)$ is any path in $D$
beginning at $x_0$ and ending at $x$.
\\

\begin{defi}\label{defturningpoint}
The zeros of the function $$g_+(x)g_-(x)=\frac{m^2c^4-(V(x)-E)^2}{c^2},$$
are called the turning points of the system (\ref{eqDirac2}).
\end{defi}
\begin{defi}\label{defstokesline}
For $x$ fixed in $D$, the set
$$\left\{y\in D,\;\re\int_x^y(g_+(t)g_-(t))^\frac{1}{2}dt=0\right\}$$
is called the Stokes line passing through $x$.
\end{defi}
We look for solutions of the form $e^{\pm\frac zh}{\Ti w}_{\pm}(z)$.
We note that due to the possible presence of such turning points,
the square root in the definition of $z(x)$ might be defined only
locally.  By formal calculations, the amplitude vector ${\Ti
w}_{\pm}(z)$ has to satisfy
\begin{eqnarray}
\frac hi\dz {\Ti w}_{\pm}(z)=\left(\begin{array}{cc}
\pm i&H(z)^{-2}\\
-H(z)^{2}&\pm i
\end{array}\right){\Ti w}_{\pm}(z).
\end{eqnarray}
The function $H(z(x))$ is given by
\begin{eqnarray}\label{H(z(x))}
H(z(x))=\left(\frac{g_-(x)}{g_+(x)}\right)^{1/4}=\left(\frac{-mc^2+(V(x)-E)}{-mc^2-(V(x)-E)}\right)^{1/4},
\end{eqnarray}
for $z(x)$ in an open simply-connected domain of the $z$-plane, where $z\longrightarrow H(z)$ is well-defined and analytic.

In order to obtain a decomposition with respect to image and kernel of the previous system
, we conjugate by
$$P_{\pm}(z)=\frac12\left(\begin{array}{cc}
H(z)&\pm iH(z)^{-1}\\
H(z)&\mp iH(z)^{-1}
\end{array}\right),\;\;\;P_{\pm}^{-1}(z)=\left(\begin{array}{cc}
 H(z)^{-1}&H(z)^{-1}\\
\mp iH(z)&\pm iH(z)
\end{array}\right),$$
and obtain a system for $w_{\pm}(z)=P_\pm(z){\Ti w}_\pm(z),$
\begin{eqnarray}
\dz { w}_{\pm}(z)=\left(\begin{array}{cc}
0&\frac{H'(z)}{H(z)}\\
\frac{H'(z)}{H(z)}&\mp\frac2h
\end{array}\right){ w}_{\pm}(z),
\end{eqnarray}
where $H'(z)$ is shorthand for $\frac{d}{dz}H(z)$. The series ansatz
\begin{eqnarray}\label{eqseries}
w_\pm(z)=\sum_{n\geq0}\left(\begin{array}{c}
w_{2n,\pm}(z)\\
w_{2n+1,\pm}(z)
\end{array}\right),
\end{eqnarray}
with $w_{0,\pm}=1$ and, for $n\geq1$, the recurrence equations
\begin{eqnarray}
\left(\dz\pm\frac2h\right)w_{2n+1,\pm}(z)&=&\frac{H'(z)}{H(z)}w_{2n,\pm}(z),\\
\dz w_{2n+2,\pm}(z)&=&\frac{H'(z)}{H(z)}w_{2n+1,\pm}(z),
\end{eqnarray}
give us a formal solution up to some additive constants. The solutions are fixed by setting
$$w_{n,\pm}({\Ti z})=0,\;\;n\geq 1,$$
at a base point ${\Ti z}=z({\Ti x})$ where ${\Ti x}\in D$ is not a turning point. We note that the previous equations for $w_{n,\pm}$
 are similar to the ones obtained by a complex WKB construction for scalar Schr\"odinger equations. See for example the works of
 C. G\'erard and Grigis \cite{CGAG} or Ramond \cite{TR}.

Let $\Om=\Om(E)$ be a simply connected subset of $D$ which does not
contain any turning point. Then the function $z=z(x)$ is conformal
from $\Om$ onto $z(\Om)$. Assume that ${\Ti z}\in z(\Om)$. If
$\Gamma_\pm({\Ti z},z)$ denotes a path of finite length  in $z(\Om)$
connecting ${\Ti z}$ and $z\in z(\Om)$, we can formally rewrite the
above differential equations for $n\geq0$ as
\begin{eqnarray}
w_{2n+1,\pm}(z)&=&\int_{\Gamma_\pm({\Ti z},z)}\exp(\pm\frac2h(\zeta-z))\frac{H'(\zeta)}{H(\zeta)}w_{2n,\pm}(\zeta)
d\zeta,\nonumber\\\nonumber
w_{2n+2,\pm}(z)&=&\int_{\Gamma_\pm({\Ti z},z)}\frac{H'(\zeta)}{H(\zeta)}w_{2n+1,\pm}(\zeta)
d\zeta,
\end{eqnarray}
or after iterated integrations, as
\begin{eqnarray}
w_{2n+1,\pm}(z)=\int_{\Gamma_\pm({\Ti z},z)}\int_{\Gamma_\pm({\Ti z},\zeta_{2n+1})}&\cdots&
\int_{\Gamma_\pm({\Ti z},\zeta_2)}\exp\left(\pm\frac2h(\zeta_1-\zeta_2+\cdots+\zeta_{2n+1}-z)\right)\times\nonumber\\\nonumber
&\times&\frac{H'(\zeta_1)}{H(\zeta_1)}\cdots\frac{H'(\zeta_{2n+1})}{H(\zeta_{2n+1})} \;d\zeta_1\cdots d\zeta_{2n+1},\\
w_{2n+2,\pm}(z)=\int_{\Gamma_\pm({\Ti z},z)}\int_{\Gamma_\pm({\Ti z},\zeta_{2n+2})}
&\cdots&\int_{\Gamma_\pm({\Ti z},\zeta_2)}\exp\left(\pm\frac2h(\zeta_1-\zeta_2+\cdots-\zeta_{2n+2})\right)\times\nonumber\\\nonumber
&\times&\frac{H'(\zeta_1)}{H(\zeta_1)}\cdots\frac{H'(\zeta_{2n+2})}{H(\zeta_{2n+2})} \;d\zeta_1\cdots d\zeta_{2n+2}.
\end{eqnarray}
\subsection{Convergence, $h$-dependence and Wronskians} We now give to the preceding formal construction some mathematical
meaning in simply connected, turning point-free compact sets
$\Om\subset D.$
\begin{lem}
For any fixed $h>0$, the series (\ref{eqseries}) converges uniformly in any compact subset of $\Om$, and
\begin{eqnarray}\label{weven}
w_\pm^{\emph{even}}(x,h)=\sum_{n\geq0}w_{2n,\pm}(z(x)),\;\;\;w_\pm^{\emph{odd}}(x,h)=\sum_{n\geq0}w_{2n+1,\pm}(z(x)),
\end{eqnarray}
are holomorphic functions in $\Om$.
\end{lem}
\begin{proof}
By assumption on $\Om$ and on $V$, the functions $w_{n,\pm}$ are
well-defined analytic functions in $\Om$. For compact subsets
$K\subset\Om$  and ${\Ti z},z\in\,z(K)$ there exist positive
constants $C_\pm^h(K)>0,$ depending on the semi-classical parameter
$h$ and the compact $K$ such that
\begin{eqnarray}
\sup_{\zeta\in\Gamma_\pm({\Ti z},z)}\Big|\exp(\pm\frac2h\zeta)\frac{H'(\zeta)}{H(\zeta)}\Big|\leq C_\pm^h(K).
\end{eqnarray}
If we denote the maximal length of the paths $\Gamma_\pm({\Ti z},\cdot)\subset K$ in the preceding iterated integrations by $$
L=\max_{{\Ti z},z\in\,z(K)}\min_{\gamma({\Ti z},z)}|\gamma({\Ti z},z)| <\infty,$$ then
$$\sup_{z\in z(K)}|w_{n,\pm}(z)|\leq\frac{C_\pm^h(K)^nL^n}{n!},\;\;\;n\geq0,$$
where the bound $\frac{L^n}{n!}$ comes from the volume of a simplex
with length $L$. Then, the lemma follows.
\end{proof}
Thus, we have  the uniform convergence of the series
(\ref{eqseries}) for $w_\pm(z)$ and complex solutions
\begin{eqnarray}\label{eqsolution}
u_\pm(x)=e^{\pm\frac{z(x)}{h}}T_\pm(z(x))\left(\begin{array}{c}
w_\pm^\ev(x)\\
w_\pm^\od(x)
\end{array}\right),
\end{eqnarray}
of the original problem (\ref{eqDirac1}) on any turning point-free
set $\Om$, where
\begin{eqnarray}\label{eqT}
T_\pm(z)=MP_\pm^{-1}(z)&=&\left(\begin{array}{cc}
0&1\\
1&0
\end{array}\right)\left(\begin{array}{cc}
 H(z)^{-1}&H(z)^{-1}\\
\mp iH(z)&\pm iH(z)
\end{array}\right)\nonumber\\&=&\left(\begin{array}{cc}
\mp iH(z)&\pm iH(z)\\
H(z)^{-1}&H(z)^{-1}
\end{array}\right),\;\;\;\;\;\;\;z\in\,z(\Om).
\end{eqnarray}
We write these solutions $u_\pm(x)$ as
\begin{eqnarray}\label{eqsolution1}
u_\pm(x;x_0,{\Ti x}),
\end{eqnarray}
indicating the particular choice of the phase base point $x_0,$ in
(\ref{eqchange}), which defines the phase function $z(x)=z(x;x_0)$,
and the choice of the amplitude base point ${\Ti z}=z({\Ti x})$,
which is the initial point of the path
 $\Gamma_\pm({\Ti z},\cdot).$\\

\begin{defi}\label{omega+-}
For ${\Ti x}\in\,\Om$ fixed, we define $\Om_\pm=\Om_\pm({\Ti x})$ the set of all $x\in\,\Om$ such that there exists a path
$\Gamma_\pm(z({\Ti x}),z(x))$ along which $x\To\pm \emph{Re } z(x)$ increases strictly.
\end{defi}
\begin{prop}\label{propexpansion}
The functions $ w_{n,\pm}$ 
 are classical analytic symbols of order \emph{[$\frac{n+1}{2}$]} in $\Om_\pm$. The functions $w_\pm^\ev(x,h)$ and $w_\pm^\od(x,h)$
 given by the identities (\ref{weven}) are classical analytic symbols of order $0$ and $1$ respectively in $\Om_\pm$.  Moreover,
we have for any 
$N\in\N$,
\begin{eqnarray*}
w_\pm^\ev(x,h)-\sum_{n=0}^Nw_{2n,\pm}(z(x))
&=&O(h^{N+1}),\\
w_\pm^\od(x,h)-\sum_{n=0}^Nw_{2n+1,\pm}(z(x))
&=&O(h^{N+2}),
\end{eqnarray*}
uniformly in any compact subsets of $\Om_\pm$. 
 In particular,
$$w_\pm^\ev(x,h)=1+h\phi(h),\;\;\;\;w_\pm^\od(x,h)=h\phi(h).$$
Here and in all this paper, $\phi(h)$ is a classical analytic symbol
of non-negative order not necessarily the same in each expression.

\end{prop}
The proof is just the same as that of \cite[Prop. 1.2]{CGAG} and \cite[Prop. 3.3]{SFCLLN}. The key point is the following:
since the iterated integrations defining $w_{n,\pm}(z)$ contain terms of the form $\exp(\pm\frac{\zeta}{h}),$ one has to
  make sure that 
$x\longmapsto\pm\re(z(x))$  is a strictly increasing function along
the path $\Gamma_\pm({\Ti z},z)$. In other words,  the paths
$\Gamma_\pm(z({\Ti x}),z(x))$ have to intersect the Stokes lines,
that is the level curves of $x\longmapsto\re(z(x))$, transversally
in a suitable direction.
\begin{defi}\label{defwronskian}
One defines the Wronskian of two $\C^2$-valued functions
$u=(u_1,u_2),\,v=(v_1,v_2)$ as
$$\W(u,v)=u_1v_2-u_2v_1.$$
\end{defi}
\begin{rem}\label{rem wronsk ind x}
For two solutions $u$ and $v$ of the equation (\ref{eq valeu propre}),
the Wronskian $\W(u,v)$ doesn't depend on $x$ and  is zero if and
only if $u$ and $v$ are proportional.
\end{rem}
If $w=\al u+\be v$ with $\al,\,\be\in\,\C$, then
$$\al=\frac{\W(w,v)}{\W(u,v)},\;\;\;\;\be=-\frac{\W(w,u)}{\W(u,v)}.$$

Elementary computations give the following complex Wronskian
formulas for complex WKB solutions with different phase  and
amplitude base points in terms of $w_\pm^\ev$ and $w_\pm^\od$.
\begin{lem}\label{lemWronsk}
Let $x_0$ and $y_0$ be two points in $\Om=\Om(E)$. If, for given
${\Ti x}$ and ${\Ti y}$, the canonical  sets $\Om_\pm({\Ti x})$ and
$\Om_\pm({\Ti y})$ have a non-empty intersection, then for any $x\in
\Om_\pm({\Ti x})\cap\Om_\pm({\Ti y})$ one has
\begin{eqnarray}
\W(u_\pm(x;x_0,{\Ti x})\!\!\!\!\!\!&,&\!\!\!\!\!\!u_\pm(x;y_0,{\Ti y}))=\pm2i\exp\left(\pm\frac1h(z(x;x_0)+z(x;y_0))\right)\\
&\times&\left(w_\pm^{\emph{odd}}(x;x_0,{\Ti
x})w_\pm^{\emph{even}}(x;y_0,{\Ti y})-
w_\pm^{\emph{even}}(x;x_0,{\Ti x})w_\pm^{\emph{odd}}(x;y_0,{\Ti
y})\right).\nonumber
\end{eqnarray}
If, for given ${\Ti x}$ and ${\Ti y}$ the canonical sets
$\Om_\pm({\Ti x})$ and $\Om_\mp({\Ti y})$ have  a non-empty
intersection, then for any $x\in \Om_\pm({\Ti x})\cap\Om_\mp({\Ti
y})$ one has
\begin{eqnarray}
\W(u_\pm(x;x_0,{\Ti x})\!\!\!\!\!\!&,&\!\!\!\!\!\!u_\mp(x;y_0,{\Ti y}))=\pm2i\exp\left(\pm\frac1h(z(x;x_0)-z(x;y_0))\right)\\
&\times&\left(w_\pm^{\emph{odd}}(x;x_0,{\Ti
x})w_\mp^{\emph{odd}}(x;y_0,{\Ti y})- w_\pm^{\emph{even}}(x;x_0,{\Ti
x})w_\mp^{\emph{even}}(x;y_0,{\Ti y})\right).\nonumber
\end{eqnarray}
\end{lem}

\section{Jost solutions}\label{sectionJost}
The Jost solutions of $Hu=Eu,$ are characterized by the behavior of the solutions at infinity.
We construct here the Jost solutions copying the procedure described
in Section \ref{sectionconst},  the new point here being that the
solutions we seek are normalized at infinity. In all this section we
will work in two unbounded, simply-connected domains $\Om^-(E)$,
$\Om^+(E)$, where $\re(V(x)+mc^2)<E$, $\re(V(x)-mc^2)>E$ respectively 
 and
 which coincide with $\calS$ for $\re x$ sufficiently large.
 The existence of such domains is of course an easy consequence of the behavior of $V$
 at infinity in $\calS$ (see assumption (A)).


First we define the phase functions with base point at infinity,
\begin{eqnarray}\label{eqchangeinfini}
z(x,\pm\infty)=\int_{\pm\infty}^x\left(\frac{m^2c^4-(V(t)-E)^2}{c^2}\right)^{1/2}\!\!\!\!\!&-&\!\!\!
\left(\frac{m^2c^4-(V^\pm-E)^2}{c^2}\right)^{1/2}dt\\\nonumber&+&\left(\frac{m^2c^4-(V^\pm-E)^2}{c^2}\right)^{1/2}x.
\end{eqnarray}
 We also see that the integral converges absolutely, hence
$$z(x,\pm\infty)=\left(\frac{m^2c^4-(V^\pm-E)^2}{c^2}\right)^{1/2}x+o(1),\;\;\;\;\;(x\To\pm\infty).$$
If the determination of the square root in $z(\cdot,\cdot)$ are the same, we get the following equalities
\begin{equation}
\label{formulephase}
\begin{aligned}
z(t_1,\pm\infty)=&z(x,\pm\infty)-z(x,t_1)=z(t_2,\pm\infty)-z(t_2,t_1)\\
z(t_1,-\infty)-z(t_1,+\infty)=&z(x,-\infty)-z(x,+\infty)=z(t_2,-\infty)-z(t_2,+\infty),
\end{aligned}
\end{equation}
where $z(\cdot,\cdot)$ is defined in (\ref{eqchange}), (\ref{eqchangeinfini}) and $x,\ t_1,\ t_2\in D$.\\

Next we define the amplitudes based at infinity. We will only define the amplitudes at $+\infty$ since the situation is similar at $-\infty$. 
 As in
Section 3 of \cite{TR}, we choose infinite paths $\gamma_\pm(x)$
starting from infinity and ending at $x$, which are asymptotically
like lines of the form $\{\mbox{Im}x=\mp\rho\mbox{Re}x\}$ for some
$\rho>0$, such that
 $x\longmapsto\mp\mbox{Re}z(x)$ are strictly increasing functions along $\gamma_\pm(x)$.
  Denoting the path $z(\gamma_\pm(x))$ by $\Gamma_\pm(+
\infty,z(x))$ and setting $w_{0,\pm}\equiv1$, we inductively define $w_{n,\pm}(z)$ by
\begin{eqnarray}
w_{2n+1,\pm}(z)&=&\int_{\Gamma_\pm(+\infty,z)}\exp(\pm\frac2h(\zeta-z))\frac{H'(\zeta)}{H(\zeta)}w_{2n,\pm}(\zeta)
d\zeta,\nonumber\\\nonumber
w_{2n+2,\pm}(z)&=&\int_{\Gamma_\pm(+\infty,z)}\frac{H'(\zeta)}{H(\zeta)}w_{2n+1,\pm}(\zeta)
d\zeta,\;\;\;n\geq0.
\end{eqnarray}
Noticing that
$$\frac{H'(x)}{H(x)}=\frac{mc^2}{2}\frac{V'(x)}{(V(x)-E)^2-m^2c^4}=O(\langle x\rangle^{-\de}),\;\;\de>1,\;\mbox{as}\;|x|\To\infty,$$
one constructs well-defined complex WKB solutions $u_{r,l}^\pm$ corresponding to these base points,
proceeding as in Section \ref{sectionconst}. 
 Here, $l$ and $r$ stand for $left$ and $right$ and correspond respectively to $x\to-\infty$ and $x\to+\infty$.
Up to a constant pre-factor, $u_{r,l}^\pm(x)$ are the previously
defined Jost solutions:
\begin{lem}\label{lemsolutionjost}
Let $u_{r,l}^\pm(x)$ be the complex WKB solutions with phase and amplitude base point at infinity. Then 
\begin{eqnarray}\label{soltuionjostdroite}
u_{r}^\pm(x)&\thicksim&\exp(\pm\frac{1}{hc}(m^2c^4-(V^+-E)^2)^{1/2}x)\left(\begin{array}{c}
\mp i\al^+\\
1/\al^+
\end{array}\right),\;\;\;x\To+\infty,\\
u_{l}^\pm(x)&\thicksim&\exp(\pm\frac{1}{hc}(m^2c^4-(V^--E)^2)^{1/2}x)\left(\begin{array}{c}
\mp i\al^-\\
1/\al^-
\end{array}\right),\;\;\;x\To-\infty,\label{soltuionjostgauche}
\end{eqnarray}
with$$\al^{\pm}=\left(\frac{V^{\pm}-E-mc^2}{-V^{\pm}+E-mc^2}\right)^{1/4}.$$
\end{lem}
\begin{proof} We just check the asymptotic behavior of $u_{r,l}^\pm(x)$
at infinity. Since $H(z(x))\To\al^{\pm}$ as $x\To\pm\infty$,  using
(\ref{eqsolution}) and (\ref{eqT}) we get by an elementary
calculation
$$u_{r}^\pm(x)\thicksim\exp(\pm\frac{1}{hc}(m^2c^4-(V^{+}-E)^2)^{1/2}x)\left(\begin{array}{cc}
\mp i\al^{+}&\pm i\al^{+}\\
1/\al^{+}&1/\al^{+}
\end{array}\right)\left(\begin{array}{c}
1\\
0
\end{array}\right),\ x\To+\infty,$$
$$u_{l}^\pm(x)\thicksim\exp(\pm\frac{1}{hc}(m^2c^4-(V^{-}-E)^2)^{1/2}x)\left(\begin{array}{cc}
\mp i\al^{-}&\pm i\al^{-}\\
1/\al^{-}&1/\al^{-}
\end{array}\right)\left(\begin{array}{c}
1\\
0
\end{array}\right),\ x\To-\infty.$$
This ends the proof of lemma.
\end{proof}

Let us now choose the determinations of
$(m^2c^4-(V^{\pm}-E)^2)^{1/2}$ and $\al^\pm$ according to the
intervals on the $E$-axis. This fixes the choice of $u_{l,r}^\pm$
and we can construct $\omin^\pm$, $\omout^\pm$ satisfying
(\ref{ominpmintro}) and (\ref{omoutintro}).
\\
1. For $E\in$ I, we choose $(m^2c^4-(V^{\pm}-E)^2)^{1/2}\in i
\R^+$, $\al^\pm\in e^{i\pi/4}\R^+$ and we denote
\begin{eqnarray}\label{ominomoutI}\omin^-:=e^{i\pi/4}u_l^+,\ \ \ \omin^+:=-e^{i\pi/4}u_r^-,\ \ \
\omout^+:=e^{i\pi/4}u_r^+,\ \ \ \omout^-:=-e^{i\pi/4}u_l^-.
\end{eqnarray}
 2. For
$E\in$ II, we choose $(m^2c^4-(V^{-}-E)^2)^{1/2}\in i \R^+$,
$\al^-\in e^{i\pi/4}\R^+$, $(m^2c^4-(V^{+}-E)^2)^{1/2}\in \R^-$,
$\al^+\in\R^+$ and we denote
\begin{eqnarray}\label{ominomoutII}\omin^-:=e^{i\pi/4}u_l^+,\ \ \
\omout^-:=-e^{i\pi/4}u_l^-,\ \ \ \om_d^+:=u_r^+.
\end{eqnarray}
 3. For $E\in$ III,
we choose $(m^2c^4-(V^{\pm}-E)^2)^{1/2}\in i \R^+$, $\al^\pm\in
e^{\mp i\pi/4}\R^+$ and we denote
\begin{eqnarray}\label{ominomoutIII}\omin^-:=e^{i\pi/4}u_l^+,\ \ \ \omin^+:=-e^{-i\pi/4}u_r^+,\ \ \
\omout^+:=e^{-i\pi/4}u_r^-,\ \ \ \omout^-:=-e^{i\pi/4}u_l^-.
\end{eqnarray}
 4. For
$E\in$ IV, we choose $(m^2c^4-(V^{+}-E)^2)^{1/2}\in i \R^+$,
$\al^+\in e^{- i\pi/4}\R^+$, $(m^2c^4-(V^{-}-E)^2)^{1/2}\in \R^-$,
$\al^+\in\R^+$ and we denote
\begin{eqnarray}\label{ominomoutIV}\omin^+:=-e^{-i\pi/4}u_r^+,\ \ \
\omout^+:=e^{-i\pi/4}u_r^-,\ \ \ \om_d^-:=u_l^-.
\end{eqnarray}
 5. For $E\in$ V, we choose $(m^2c^4-(V^{\pm}-E)^2)^{1/2}\in i
\R^+$, $\al^\pm\in e^{-i\pi/4}\R^+$ and we denote
\begin{eqnarray}\label{ominomoutV}\omin^-:=e^{-i\pi/4}u_l^-,\ \ \
\omin^+:=-e^{-i\pi/4}u_r^+,\ \ \ \omout^+:=e^{-i\pi/4}u_r^-,\ \ \
\omout^-:=-e^{-i\pi/4}u_l^+.
\end{eqnarray}

\begin{thm}\label{solutioJost}
For real $E$, (\ref{eq valeu propre}) has solutions of the following form:

1. For $E\in$ \emph{I, III} or  \emph{V}, there are four Jost solutions $\omin^\pm$, $\omout^\pm$ which behave like 
\begin{eqnarray}
 \omin^\pm&\sim&\exp\{\mp\frac{i}{hc}\Phi(E-V^\pm)
x\}\left(\begin{array}{c}
A(E-V^\pm)\\
{\mp}{A(E-V^\pm)^{-1}}
\end{array}\right)\,\;\;\;\mbox{as }x\To\pm\infty,\label{omin}\\
 \omout^\pm&\sim&\exp\{\pm\frac{i}{hc}\Phi(E-V^\pm)
x\}\left(\begin{array}{c}
A(E-V^\pm)\\
{\pm}{A(E-V^\pm)^{-1}}
\end{array}\right)\,\;\;\;\mbox{as }x\To\pm\infty,\label{omout}
\end{eqnarray}
with $\Phi(E)=\emph{sgn}(E)\sqrt{E^2-m^2c^4}$, $A(E)=\sqrt[4]{\frac{E+mc^2}{E-mc^2}}$ and $\emph{sgn}(E)=\frac{E}{|E|}$
for $E\not\in[-mc^2,mc^2]$. 

2. For $E\in$ \emph{II} (resp. $E\in$ \emph{IV}), there are two Jost
solutions $\omin^-$, $\omout^-$ (resp. $\omin^+$, $\omout^+$) which
behave as in (\ref{omin}), (\ref{omout})  and a decreasing solution
$\om_d^+$ (resp. $\om_d^-$). The solutions $\om_d^\pm$ behave
exactly like
\begin{eqnarray}
 \om_d^\pm&\sim&\exp\{\mp\frac{1}{hc}\sqrt{m^2c^4-(V^\pm-E)^2}x\}\left(\begin{array}{c}
\mp i\sqrt[4]{\frac{mc^2+E-V^\pm}{mc^2-E+V^\pm}}\\
\sqrt[4]{\frac{mc^2-E+V^\pm}{mc^2+E-V^\pm}}
\end{array}\right)\,\;\;\;\mbox{as }x\To\pm\infty.
\end{eqnarray}

%
%

\end{thm}

For $E\in$ I, III or V, according to the relation (\ref{matrixS}), it is
sufficient to calculate the two terms $r(E,h),\ t(E,h)$ in $\bbT$ to
obtain the matrix $\bbS$.
 The definition of the Wronskian (see Definition \ref{defwronskian}) leads to:
\begin{eqnarray}\label{eq1tr}
t(E,h)&=&\frac{\W(\omin^-,\omin^+)}{\W(\omout^+,\omin^+)},\\
r(E,h)&=&\frac{\W(\omout^-,\omin^+)}{\W(\omout^+,\omin^+)}.\label{eq2tr}
\end{eqnarray}
\section{The  Klein paradox case}\label{sect-2ptourn}
We suppose that $V$ satisfies assumption (A), the energy $E\in$ III
and $m>0$ (see Fig.  2). In this section we will work in two
unbounded, simply-connected domains $\Om^-(E)$, $\Om^+(E)$, where
$\re(V(x)+mc^2)<E$, $\re(V(x)-mc^2)>E$ respectively and
 which coincide with $\calS$ for $|\re x|$ sufficiently large.
 Using Theorem \ref{solutioJost}, Proposition \ref{propexpansion} and (\ref{ominomoutIII}) there are two Jost solutions in
 $\Om^\pm(E)$:
 \begin{eqnarray}\label{ominpm}
 \omin^\pm&=&\exp\{\frac1hz(x,\pm\infty)\}\left(\begin{array}{c}
\Ti H(z(x))\\
{\mp}{\Ti H(z(x))^{-1}}
\end{array}\right)(1+h\phi(h))
\\
\label{omoutpm}
 \omout^\pm&=&\exp\{\frac{-1}{h}z(x,\pm\infty)\}\left(\begin{array}{c}
\Ti H(z(x))\\
{\pm}{\Ti H(z(x))^{-1}}
\end{array}\right)(1+h\phi(h)).
 \end{eqnarray}
 The function $z(x,\pm\infty)$ is defined by
(\ref{eqchangeinfini}) and \begin{eqnarray}\label{TiH(z(x))}\Ti
H(z(x))=\left(\frac{E-V(x)+mc^2}{E-V(x)-mc^2}\right)^\frac14.
\end{eqnarray} On $\Om^\pm(E)\cap\R$, we have:
\begin{eqnarray}\label{z(x,infty determination+)}z(x,\pm\infty)&=& i\int_{\pm\infty}^{x}\sqrt{\frac{(E-V(t))^2-m^2c^4}{c^2}}-\sqrt{\frac{(E-V^\pm)^2-m^2c^4}{c^2}}\ dt\\
& &+i\sqrt{\frac{(E-V^\pm)^2-m^2c^4}{c^2}}\ x\ \ \ \ \ \ \
\nonumber\\\label{TiH(z(x))real} \Ti
H(z(x))&=&\sqrt[4]{\frac{E-V(x)+mc^2}{E-V(x)-mc^2}}.
\end{eqnarray}

We suppose that there are only two real turning points
$t_1(E)<\,t_2(E)$ and that they are simple. Notice that $t_1(E)$ is
a zero of $E-V(t)-mc^2$ and $t_2(E)$ is a zero of $E-V(t)+mc^2$. In
that case the Stokes lines are as shown in Fig. 4. In order to
obtain $\bbS$, we compute the Wronskians given in (\ref{eq1tr}),
(\ref{eq2tr}) and then the coefficients $t(E,h)$, $r(E,h)$.
\\
{\begin{center}
\psfrag{x2-}[rh][][1][0]{}
\psfrag{t2}[][][1][0]{$\ \ \ t_2(E)$}
\psfrag{x2+}[][][1][0]{}
\psfrag{t1}[][][1][0]{$\!\!\!t_1(E)$}
\psfrag{x1-}[][][1][0]{}
\psfrag{x1+}[][][1][0]{}
\psfrag{ga1}[][][1][0]{$\gamma_1$}
\psfrag{ga2}[][][1][0]{$\gamma_2$} \psfrag{titre}[][][1][0]{Fig. 4.
The turning points and the paths $\gamma_j$}
\includegraphics[height=5cm, width=14cm]
{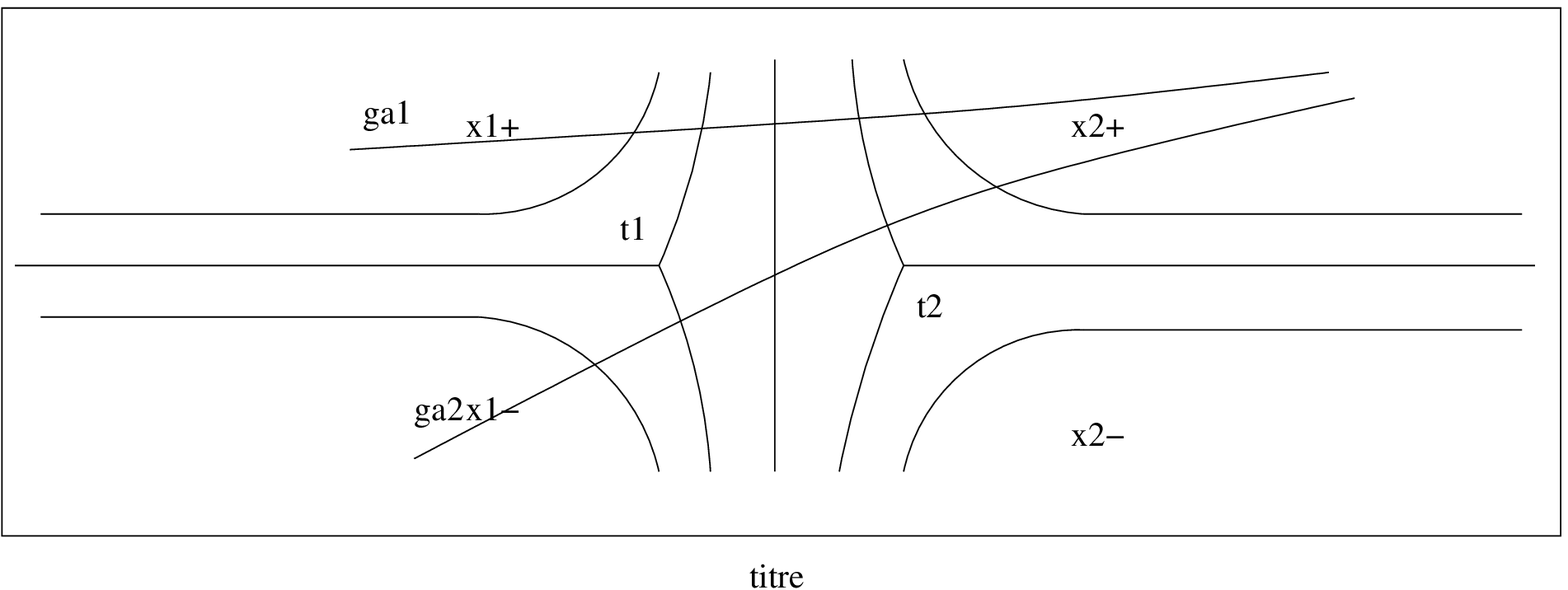}
\end{center}
}
\vspace{3mm}

{\medskip\it Computation of $\W(\omout^+,\omin^+)$:} Since the two
solutions $\omout^+,\omin^+$ are defined in $\Om^+(E)$, we can
compute this Wronskian in $\Om^+(E)$ and from Lemma \ref{lemWronsk}
we obtain
\begin{eqnarray}\label{wronskianequal-2}
\W(\omout^+,\omin^+)=-2
.\end{eqnarray}

\medskip {\it Computation of $\W(\omin^-,\omin^+)$:} The two
solutions $\omin^-,\omin^+$ are defined in $\Om^-(E),\Om^+(E)$
respectively. Since the Wronkians $\W(\omin^-,\omin^+)(x)$ are
independent on $x$ (see Remark \ref{rem wronsk ind x}) we compute
this Wronskian in $\Om^-(E)$ for example. For that we extend
$\omin^+$ , which is defined in
 $\Om^+(E)$, into $\Om^-(E)$. We will extend the square
root in $\omin^+$ which is defined in $\Om^+(E)$, into
$\C\setminus\{-\im(z)>0,\ \re(z)=t_{1}(E)\}\cup\{-\im(z)>0,\
\re(z)=t_{2}(E)\}$. Thanks to the structure of the Stokes lines
between $t_1(E)$ and $t_2(E)$, we can find a path $\gamma_1$ from
$+\infty(1+i\de_1)$
 to $-\infty(1-i\de_1)$ (for $\de_1>0$) transverse to the Stokes lines
along which we can extend $\omin^+$. We remark that between $t_1(E)$
and $t_2(E)$ on the real axis we have $(E-V(t))^2-m^2c^4<0$. The
extension of $t\in]t_2(E),+\infty[\To
\sqrt{\frac{(E-V(t))^2-m^2c^4}{c^2}}$   coincide with
$i\sqrt{\frac{m^2c^4-(E-V(t))^2}{c^2}}$ on $]t_1(E),t_2(E)[$ and
with $-\sqrt{\frac{(E-V(t))^2-m^2c^4}{c^2}}$ on $]-\infty,t_1(E)[$.
On the other hand, the extension of $\Ti H(z(x))$ stay in $\R^+$ on
$]-\infty,t_1(E)[$.
  If we denote by $\omin^{+,1}$ the extension of
$\omin^+$ along $\gamma_1$, we have:
\begin{eqnarray*}
\omin^{+,1}&=&\exp\{\frac1h(-z(x,t_1(E))+z(t_2(E),+\infty)+S(E))\}\left(\begin{array}{c}
\Ti H(z(x))\\
{-}{\Ti H(z(x))^{-1}}
\end{array}\right)(1+h\phi(h)),
\end{eqnarray*}
 with  $z(t_2(E),+\infty)$ defined in (\ref{z(x,infty
 determination+)}) and
\begin{eqnarray}\label{z(x,t_1)}
z(x,t_1(E))&=&i\int_{t_1(E)}^x\left(\frac{(E-V(t))^2-m^2c^4}{c^2}\right)^{\frac12}dt,\\
S(E)&=&\int_{t_1(E)}^{t_2(E)}\sqrt{\frac{m^2c^4-(E-V(t))^2}{c^2}}dt.\label{S(E)}
\end{eqnarray}
Here, $\left(\frac{(E-V(t))^2-m^2c^4}{c^2}\right)^{\frac12}\in\R^+$
for $t\in]-\infty,t_1(E)[$.

Then,
\begin{eqnarray*}
\W(\omin^-,\omin^+)=-2(1+h\phi(h))\exp\{\frac1h(z(t_1(E),-\infty)+z(t_2(E),+\infty)+S(E))\},
\end{eqnarray*}
where $z(t_2(E),+\infty),\ z(t_1(E),-\infty)$ are defined in
(\ref{z(x,infty
determination+)}).

\medskip
{\it Computation of $\W(\omout^-,\omin^+)$:} This wronskian is also
between two solutions which are defined in different domains, then
we extend one of these solutions into the domain of the other
solution. For example we extend $\omin^+$, which is defined in
$\Om^+(E)$, into $\Om^-(E)$ which is a subset of
$\C\setminus\{\im(z)>0,\ \re(z)=t_{1}(E)\}\cup\{-\im(z)>0,\
\re(z)=t_{2}(E)\}$. Here, we can also find a path $\gamma_2$ from
$+\infty(1+i\de_2)$ to $-\infty(1+i\de_2)$ for $\de_2>0$ transverse
to the  Stokes lines along which we can extend $\omin^+$ into
$\Om^-(E)$. If we denote by $\omin^{+,2}$ the extension of $\omin^+$
along $\gamma_2$, we have:
$$\omin^{+,2}=\exp\{\frac1h(z(x,t_1(E))+z(t_2(E),+\infty)+S(E))\}\left(\begin{array}{c}
i\Ti H(z(x))\\
{i}{\Ti H(z(x))^{-1}}
\end{array}\right)(1+h\phi(h)).$$
Here $\Ti H(z(x))\in\R^+$ on $]-\infty,t_1(E)[$ and $z(x,t_1(E)),\  S(E)$
are defined in (\ref{z(x,t_1)}), (\ref{S(E)}).

The computation of $\W(\omout^-,\omin^+)$  yields:
\begin{eqnarray*}
\W(\omout^-,\omin^+)=2i(1+h\phi(h))\exp\{\frac1h(-z(t_1(E),-\infty)+z(t_2(E),+\infty)+S(E))\}.
\end{eqnarray*}

Then, we obtain (see (\ref{eq1tr}) and (\ref{eq2tr})):
\begin{eqnarray*}
t(E,h)&=&(1+h\phi(h))\exp\{\frac1h(z(t_1(E),-\infty)+z(t_2(E),+\infty)+S(E))\},\\
r(E,h)&=&-i(1+h\phi(h))\exp\{\frac1h(-z(t_1(E),-\infty)+z(t_2(E),+\infty)+S(E))\}.
\end{eqnarray*}
Since  $\phi(h)$ is a classical analytic symbols of non-negative
order and using (\ref{matrixS}) we have:
\begin{eqnarray*}
s_{11}&=&\frac{1}{\overline{\,t\,}(E,h)}=(1+h\phi_1(h))\exp\{\frac1h(z(t_1(E),-\infty)+z(t_2(E),+\infty)-S(E))\},\\
s_{21}&=&\frac{\overline{\,r\,}(E,h)}{\overline{\,t\,}(E,h)}=(i+h\phi_2(h))\exp\{\frac2h(z(t_1(E),-\infty)\},\\
s_{12}&=&\frac{-r(E,h)}{\overline{\,t\,}(E,h)}=(i+h\phi_3(h))\exp\{\frac2h(z(t_2(E),+\infty)\}.
\end{eqnarray*}
The functions $\phi_1(h),\,\phi_2(h),\,\phi_3(h)$ are classical
analytic symbols of non-negative order. 
 This ends the proof of Theorem \ref{thm2point}.
\section{Total transmission}\label{sec-0ptourn}
We suppose that $V$ satisfies assumption (A), the energy $E\in$ I or
$E\in$ V and $m\geq0$ (see Fig.  5, Fig.  6).
{\begin{center}
\psfrag{E}[][][1][0]{E} \psfrag{t1}[rh][][1][0]{$t_1$}
\psfrag{t2}[][][1][0]{$t_2$}
\psfrag{Vl+mc2}[][][1][0]{$V^-\!\!\!+mc^2$}
\psfrag{Vl-mc2}[][][1][0]{$V^-\!\!\!-mc^2$}
\psfrag{Vr+mc2}[][][1][0]{$V^++mc^2$}
\psfrag{Vr-mc2}[][][1][0]{$V^+-mc^2$} \psfrag{titre1}[][][1][0]{Fig.
5. $V\pm mc^2$ and $E\in$I } \psfrag{titre2}[][][1][0]{Fig. 6. $V\pm
mc^2$ and $E\in$V} \psfrag{.}[][][1][0]{}
 \includegraphics[height=5cm, width=7cm]
{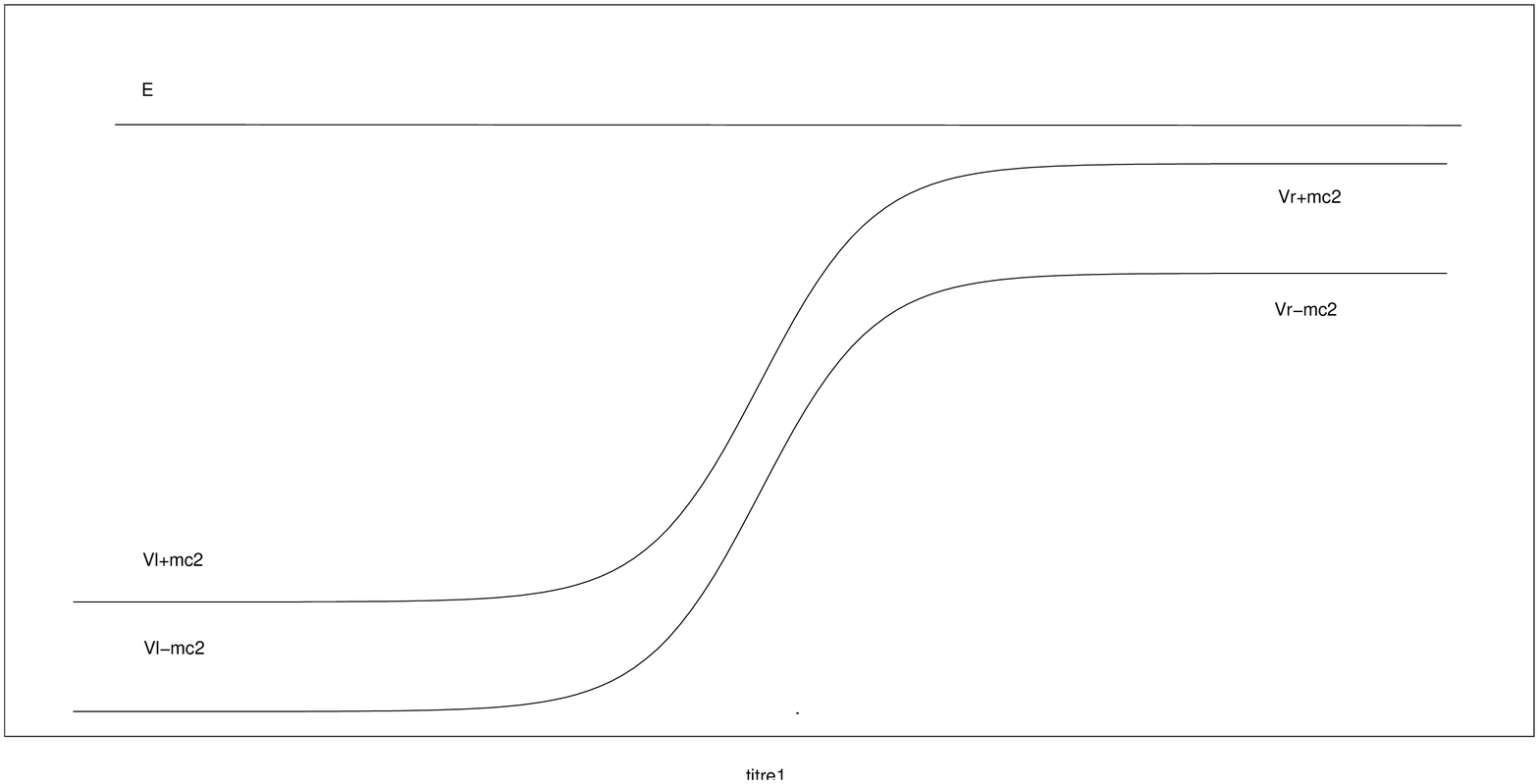}\hspace{1cm}
\includegraphics[height=5cm, width=7cm]
{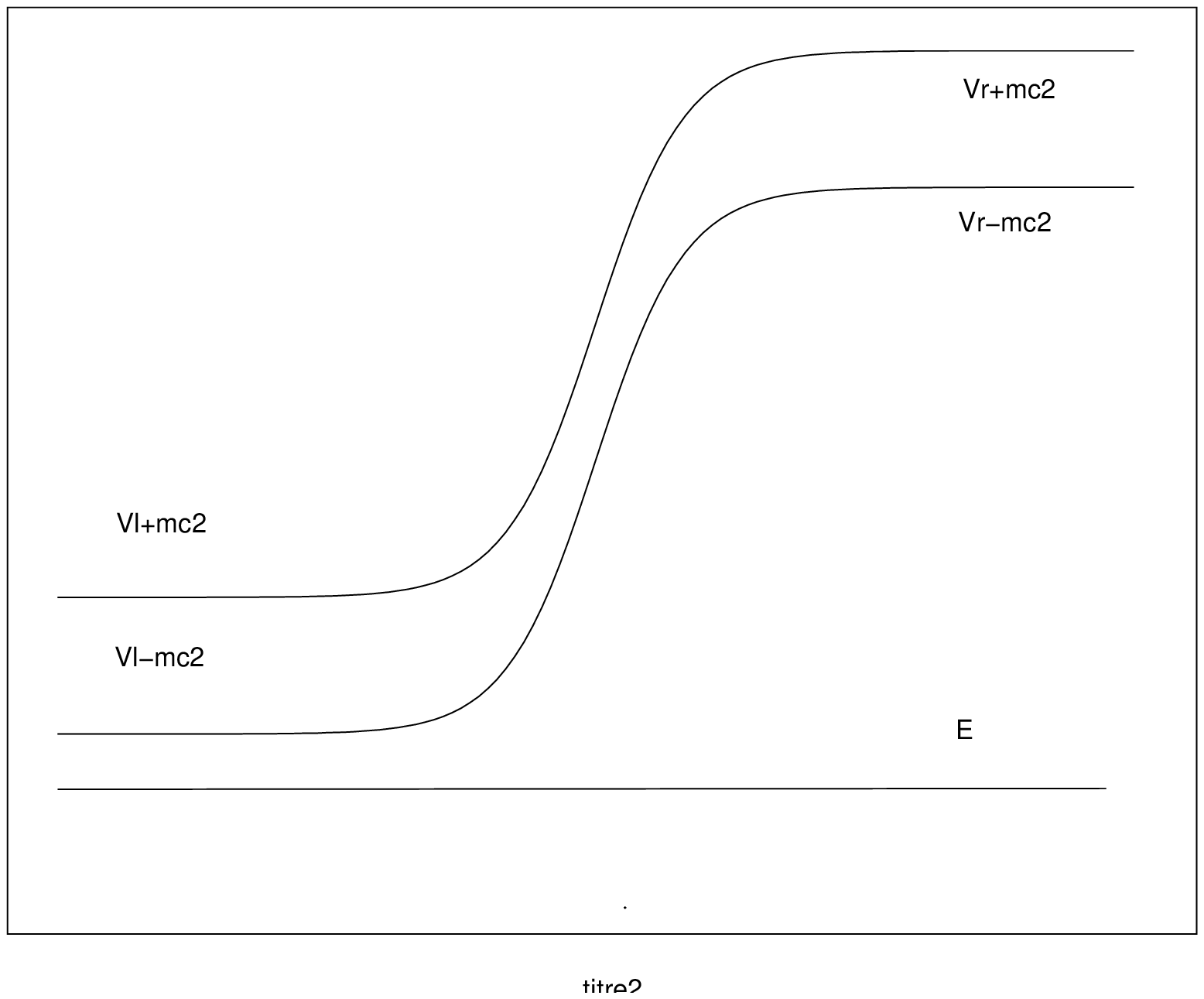}
\end{center}}
We suppose that  there exists no real turning point. In that case
the Stokes lines are horizontal lines near the real axis. We will only work for 
$E\in$ I. 
The case where $E\in$ V can be treated similarly.\\

In this section 
 we work in
$\Om^-(E)$ defined in the previous section. Now this set is a
neighborhood of the real axis. Using Theorem \ref{solutioJost},
Proposition \ref{propexpansion} and (\ref{ominomoutI}) there are
four Jost solutions:
 \begin{equation}\label{winout+-transmission}
\begin{aligned}
\omin^\pm=&\exp\{\frac{\mp1}{h}z(x,\pm\infty)\}\left(\begin{array}{c}
\Ti H(z(x))\\
{\mp}{\Ti H(z(x))^{-1}}
\end{array}\right)(1+h\phi(h))
\\
 \omout^\pm=&\exp\{\frac{\pm1}{h}z(x,\pm\infty)\}\left(\begin{array}{c}
\Ti H(z(x))\\
{\pm}{\Ti H(z(x))^{-1}}
\end{array}\right)(1+h\phi(h)).
\end{aligned}
 \end{equation}
 The functions $z(x,\pm\infty)$ and $\Ti H(z(x))$ are defined in (\ref{eqchangeinfini}) and (\ref{TiH(z(x))}) and coincide  with 
(\ref{z(x,infty determination+)}), (\ref{TiH(z(x))real}) on the real
axis. Here, the setting is different from the previous section. The
solutions $\omin^\pm$ and $\omout^\pm$ are defined in the same
domain $\Om^-(E)$ and there are no problem to extend the different
square roots.

As in Section \ref{sectionJost}, it is sufficient to calculate the
two terms $r(E,h),\ t(E,h)$ (see (\ref{eq1tr}), (\ref{eq2tr})) to
obtain the matrix $\bbS$.

\medskip
 {\it Computation of
$\W(\omin^-,\omin^+)$, $\W(\omout^+,\omin^+)$:} Since the function
$\Ti H$ in (\ref{winout+-transmission}) is the same for $\omin^-$
and $\omin^+$, we have:

\begin{eqnarray}\label{wronskianwin-wout+transmission}
\W(\omin^-,\omin^+)(x)&=&-2(1+h\phi(h))\exp\{\frac1h(z(x,+\infty)-z(x,-\infty))\}\\\nonumber
&=&-2(1+h\phi(h))\exp\{\frac1h(z(0,+\infty)-z(0,-\infty))\}.
\end{eqnarray}
Moreover, as in (\ref{wronskianequal-2})
\begin{eqnarray*}
\W(\omout^+,\omin^+)(x)&=&-2
.
\end{eqnarray*}
Then, we obtain (see (\ref{eq1tr}))
\begin{eqnarray}\label{eqvaleurt(E,h)}t(E,h)=(1+h\phi(h))\exp\{\frac1h(z(0,+\infty)-z(0,-\infty))\},
\end{eqnarray}
with \begin{eqnarray}\label{eqz(0,pminfty,totaltaransmission}
z(0,\pm\infty)=
i\int_{\pm\infty}^{0}\sqrt{\frac{(E-V(t))^2-m^2c^4}{c^2}}-\sqrt{\frac{(E-V^\pm)^2-m^2c^4}{c^2}}\
dt.\end{eqnarray}
 Since  $\phi(h)$ is a classical analytic symbol of
non-negative order and using (\ref{matrixS}) we have:
\begin{eqnarray*}
s_{11}&=&\frac{1}{\overline{\,t\,}(E,h)}=(1+h\ti\phi(h))\exp\{\frac1h(z(0,+\infty)-z(0,-\infty))\}.
\end{eqnarray*}

\medskip
{\it 
 Computation of $\W(\omout^-,\omin^+)$:} As in (\ref{wronskianwin-wout+transmission}),\\
\begin{eqnarray*}
 \W(\omout^-,\omin^+)(x)=O(h)\exp\{-\frac1h(z(x,+\infty)+z(x,-\infty))\}.
\end{eqnarray*}
Using that the square root in $z(x,+\infty)$ and $z(x,-\infty)$ have the same determination, we have
\begin{eqnarray*}
 \W(\omout^-,\omin^+)(x)=O(h)\exp\{-\frac1h(z(0,+\infty)+z(0,-\infty))\}\exp\{-\frac2h(z(x,0))\},
\end{eqnarray*}
 where $z(0,\pm\infty)\in i\R$ is defined in (\ref{eqz(0,pminfty,totaltaransmission})
 and $z(x,0)=i\int_{0}^{x}\sqrt{\frac{(E-V(t))^2-m^2c^4}{c^2}}dt$.
  Since the Wronskians are independent on $x$, we estimate the term $z(x,0)$ for $x=-iy,\ 0<y\ll1$.
  Here, we have $z(x,0)=z(-iy,0)=-iy(i\sqrt{\frac{(E-V(0))^2-m^2c^4}{c^2}})+O(y^2)=Cy+O(y^2)$ for $C>0$.
  Thereafter, $\W(\omout^-,\omin^+)=O(e^{-C/h})$ for an other $C>0$ and then
\begin{eqnarray}\label{eqvaleurrr}
r(E,h)=
O(e^{-C/h}).
\end{eqnarray}

Consequently, using (\ref{matrixS}), we have, for a positive
constant $C$,
$$s_{21}=\frac{\overline{\,r\,}(E,h)}{\overline{\,t\,}(E,h)}=
O(e^{-C/h})
,$$
$$
s_{12}=\frac{-r(E,h)}{\overline{\,t\,}(E,h)}=
O(e^{-C/h}).
$$
 This ends the proof of Theorem \ref{thmtransmission}.

\section{Total reflection}\label{sec-1ptourn}

We suppose here that $V$ satisfies assumption (A), the energy $E\in$ II or IV and $m>0$  (see Fig.  7 or Fig.  8).\\
{\begin{center}
\psfrag{E}[][][1][0]{E} \psfrag{t1}[rh][][1][0]{$t_1(E)$}
\psfrag{t2}[][][1][0]{$\ \ \ \ t_2(E)$} \psfrag{Vl+mc2}[][][1][0]{$\
V^-+mc^2$} \psfrag{Vl-mc2}[][][1][0]{$\ V^--mc^2$}
\psfrag{Vr+mc2}[][][1][0]{$V^+\!+mc^2$}
\psfrag{Vr-mc2}[][][1][0]{$V^+\!-mc^2$}
\psfrag{titre1}[][][1][0]{Fig. 7. $E\in$II}
\psfrag{titre2}[][][1][0]{Fig. 8. $E\in$IV}
 \includegraphics[height=5cm, width=7cm]
{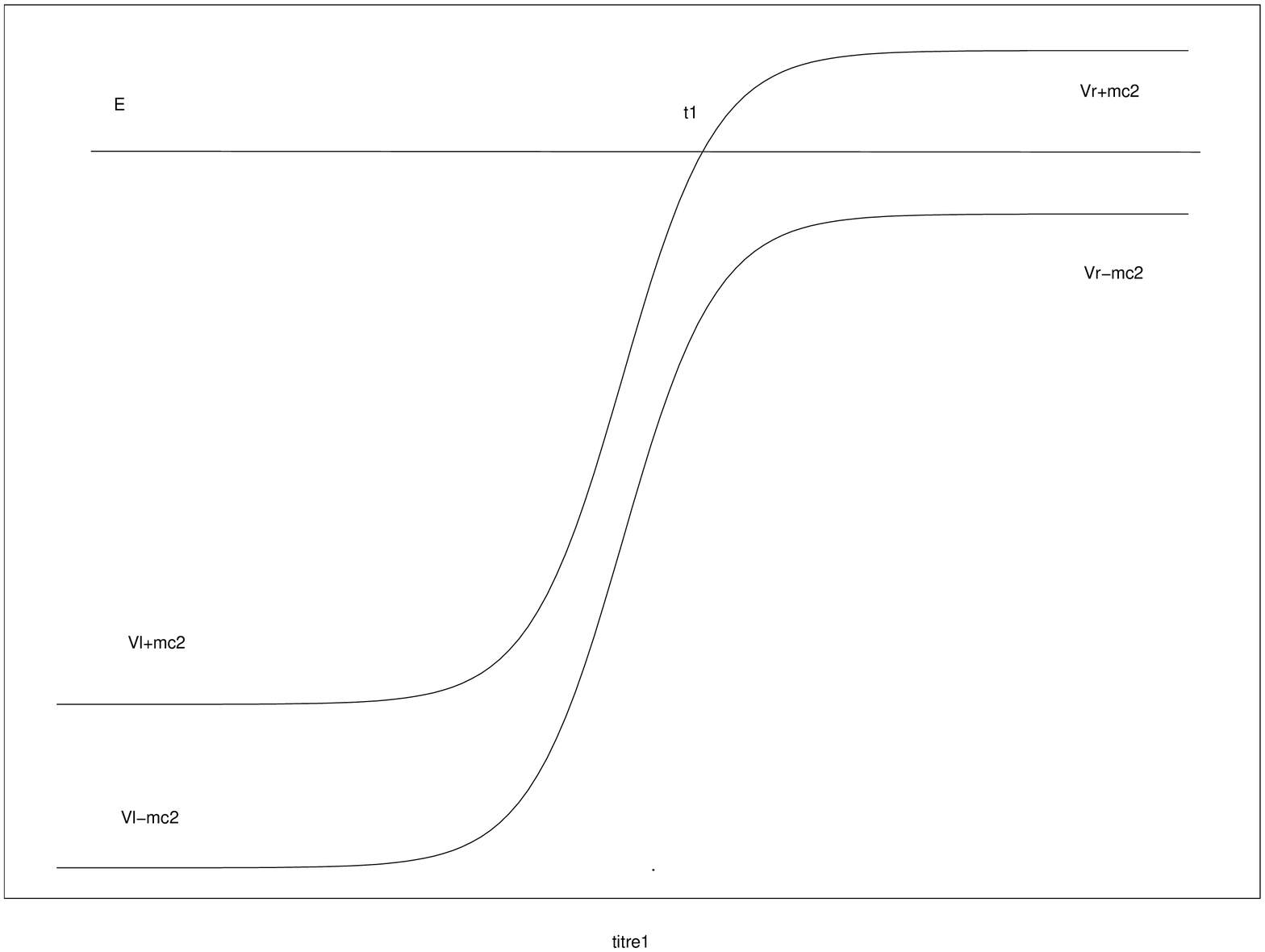}\hspace{1cm}
\includegraphics[height=5cm, width=7cm]
{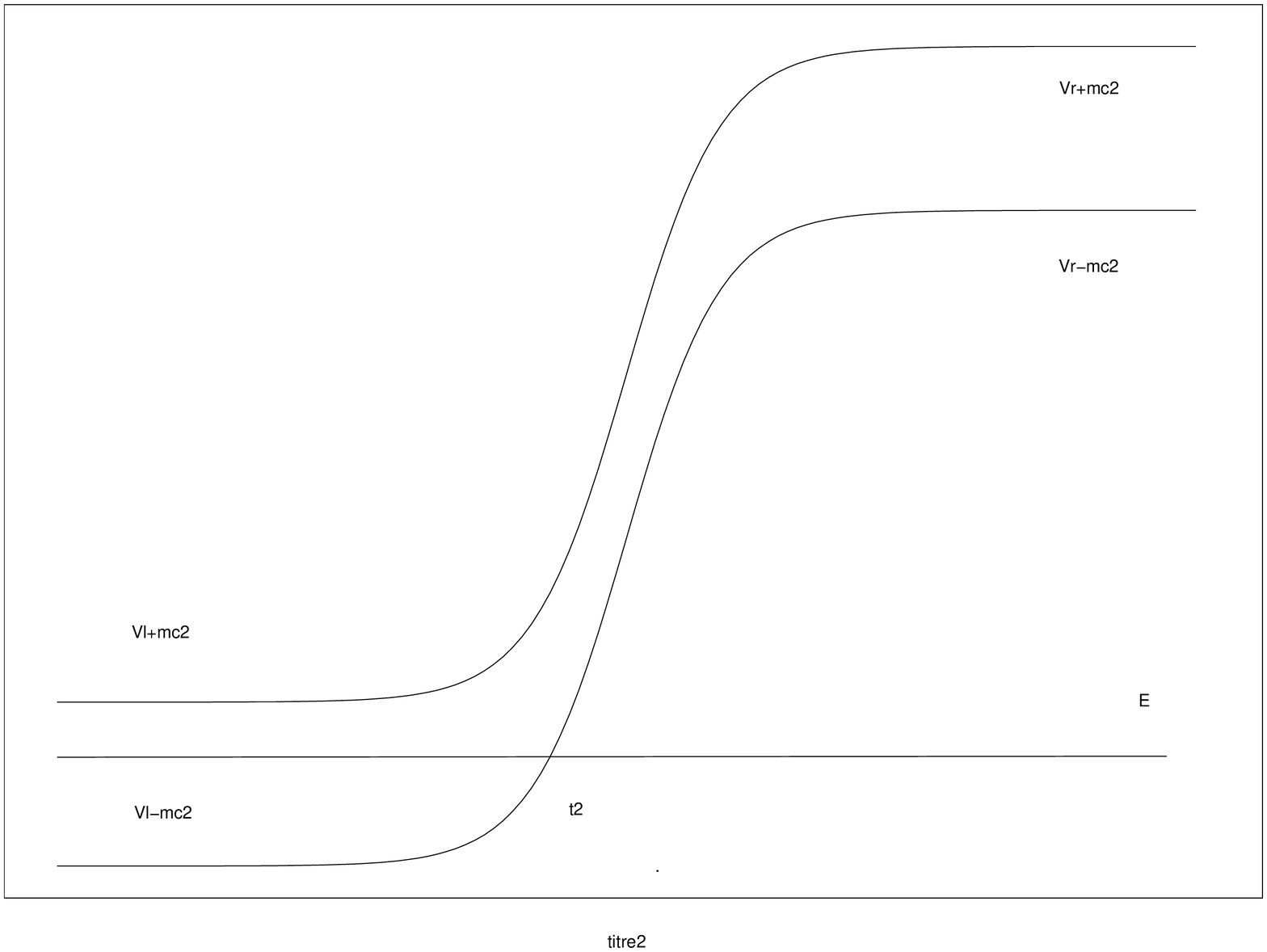}
\end{center}}
As in the Section \ref{sect-2ptourn}, we will work in two unbounded,
simply-connected domains $\Om^-(E)$, $\Om^+(E)$, where $\re
(V(x)-E+mc^2)<E$, $\re(V(x)-mc^2)>0$ respectively and which coincide
with ${\cal S}$ for $|\re x|$ sufficiently large. 
  We will only
work for $E\in$ II. The case  $E\in$ IV can be treated similarly. Using
Theorem \ref{solutioJost}, Proposition \ref{propexpansion} and
(\ref{ominomoutII}) there are two Jost solutions in
 $\Om^-(E)$:
 \begin{eqnarray*}
 \omin^-&=&\exp\{\frac1hz(x,-\infty)\}\left(\begin{array}{c}
\Ti H(z(x))\\
{+}{\Ti H(z(x))^{-1}}
\end{array}\right)(1+h\phi(h))
\\
 \omout^-&=&\exp\{\frac{-1}{h}z(x,-\infty)\}\left(\begin{array}{c}
\Ti H(z(x))\\
{-}{\Ti H(z(x))^{-1}}
\end{array}\right)(1+h\phi(h))
 \end{eqnarray*}
with $\phi(h)$ a classical analytic symbol of non-negative order.
 The functions $z(x,-\infty)$, $\Ti H(z(x))$ are defined
in (\ref{eqchangeinfini}), (\ref{TiH(z(x))}) and coincide with
(\ref{z(x,infty determination+)}), (\ref{TiH(z(x))real}) on the real
axis. 
From Lemma \ref{lemsolutionjost}, there exist an exponentially
decreasing Jost solution and an exponentially  increasing one. As
explained before Theorem \ref{thmt1pt}, we exclude the increasing
solutions which does not represent a physical state. We limit
ourself to the one-dimensional space generated by the decreasing
solution $\om_d^+$ which satisfies in $\Om^+(E)$: 
\begin{eqnarray}
\om_d^+=\exp\{\frac{-1}{h}z(x,+\infty)\}\left(\begin{array}{c}
 -iH(z(x))\\
{H(z(x))^{-1}}
\end{array}\right)(1+h\phi(h)),
\end{eqnarray}
from Theorem
\ref{solutioJost}, Proposition \ref{propexpansion} and
(\ref{ominomoutII}).
The functions $z(x,+\infty)$, $H(z(x))$ are defined in
(\ref{eqchangeinfini}), $(\ref{H(z(x))})$ respectively and coincide,
on the real axis, with
\begin{eqnarray}\label{z(x,+infini determina+ reflection)}
z(x,+\infty)&=&\frac1c\int_{+\infty}^x\sqrt{m^2c^4-(E-V(t))^2}-\sqrt{m^2c^4-(E-V^+)^2}\
dt\\
& &+\frac1c\sqrt{m^2c^4-(E-V^+)^2}\ x\nonumber\\
 H(z(x))&=&\sqrt[4]{\frac{mc^2+E-V(x)}{mc^2-E+V(x)}}.
\end{eqnarray}

We suppose that there is only one real turning point $t_1(E)$ 
and that it is simple. In that case the Stokes lines are as shown in
the fig.
Fig. 9.\\\\
{\begin{center}
\psfrag{x2-}[rh][][1][0]{}
\psfrag{t2}[][][1][0]{$t_2$}
\psfrag{x2+}[][][1][0]{}
\psfrag{t1}[][][1][0]{$t_1$}
\psfrag{x1-}[][][1][0]{}
\psfrag{x1+}[][][1][0]{}
\psfrag{titre}[][][1][0]{Fig. 9. Stokes lines for $E\in$II}
\psfrag{titre1}[][][1][0]{Fig. 10. Stokes lines for $E\in$IV}
\includegraphics[height=4cm, width=7cm]
{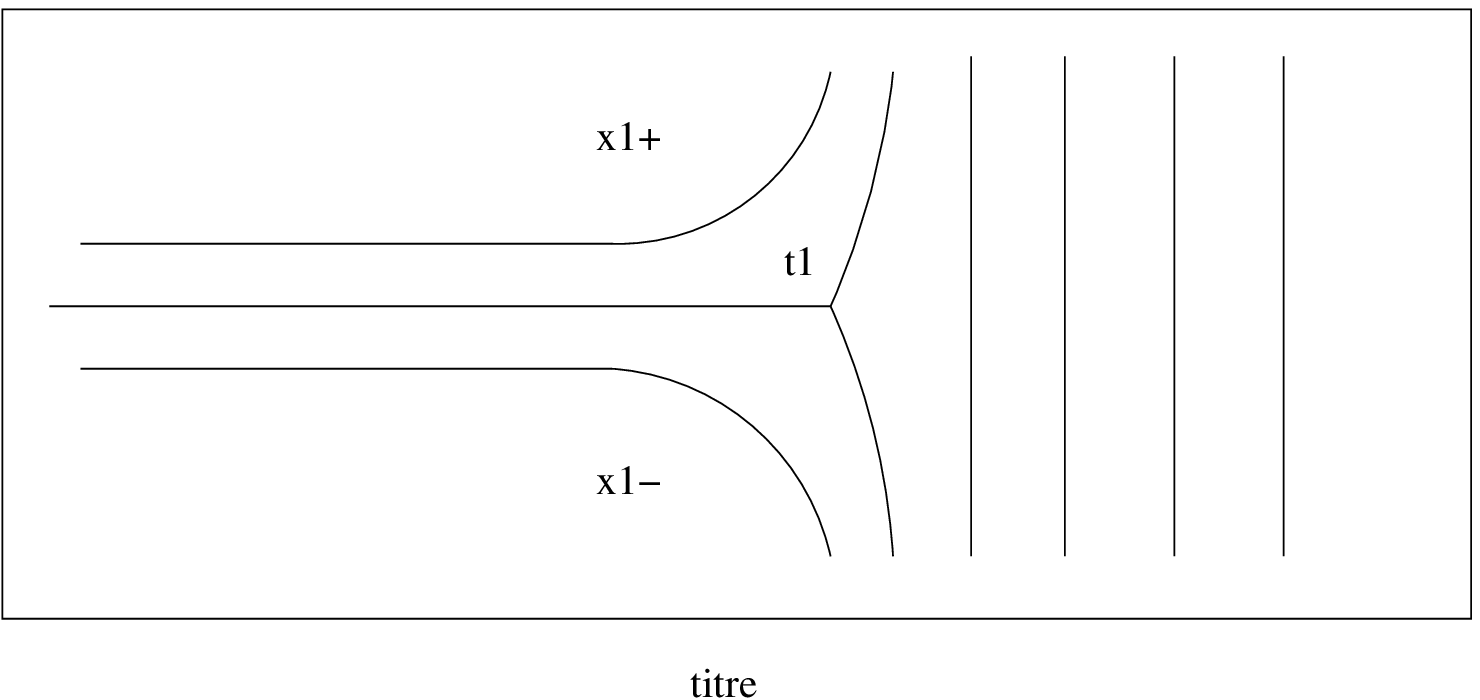}\hspace{1cm}
\includegraphics[height=4cm, width=7cm]
{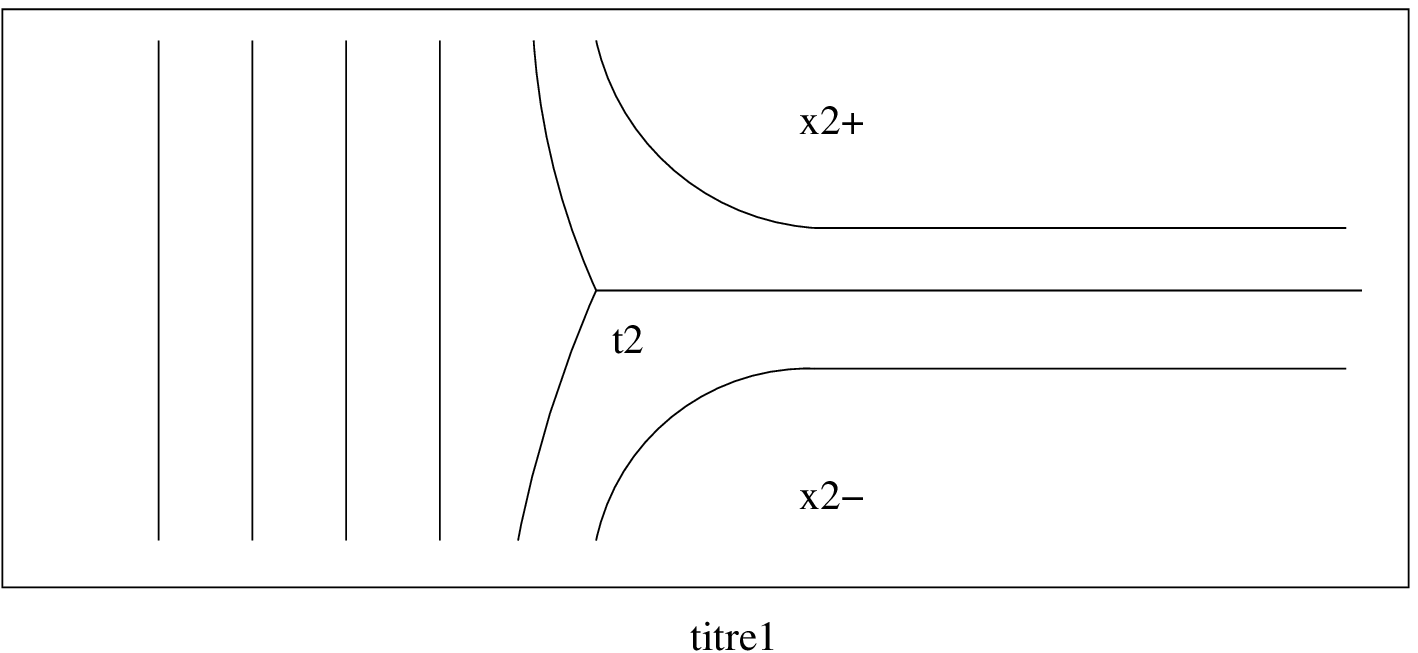}
\end{center}}

 According to the definition of
the Wronskian, 
 we have
$$\al_{\rm out}^-=\frac{\W(\omin^-,\om_d^+)}{\W(\om_d^+,\omout^-)}, \ \ \ \ \ \
\be_d^+=\frac{\W(\omin^-,\omout^-)}{\W(\om_d^+,\omout^-)}.$$

\medskip
{\it Computation of $\W(\omin^-,\om_d^+)$:} In order to calculate
this Wronskian we need to extend one of the solutions
$\omin^-,\,\om_d^+$ from its domain to the domain of the other
solution, for example, we extend $\om_d^+$ from $\Om^+(E)$ to
$\Om^-(E)$. For that, we extend the square root in $\om_d^+$ to
$\C\setminus\{-\im(z)>0,\ \re(z)=t_{1}\}$. Thanks to the structure
of the Stokes lines , we can find a path $\ti\gamma_1$ from
$+\infty(1-i\ti\de_1)$ to $-\infty(1-i\ti\de_1)$ (for $\ti\de_1>0$)
transverse to the Stokes lines along which we can extend $\om_d^-$.
The extension of $t\in]t_1(E),+\infty[\To\sqrt{m^2c^4-(E-V(t))^2}$
coincides with $i\sqrt{(E-V(t))^2-m^2c^4}$ on $]-\infty,t_1(E)[$. On
the other hand, on $]-\infty,t_1(E)[$, $H(z(x))$ takes its values in
$e^{-i\pi/4}\R^+$. If we denote by $\om_d^{+,1}$ the extension of
$\om_d^+$ along $\ti\gamma_1$, we have
$$\om_d^{+,1}=\exp\{\frac1h(-z(x,t_1(E))-z(t_1(E),+\infty))\}\left(\begin{array}{c}
 -e^{i\pi/4}\Ti H(z(x))\\
e^{i\pi/4}{\Ti H(z(x))^{-1}}
\end{array}\right)(1+h\phi(h)),$$
with $z(t_1(E),+\infty)$, $\Ti H(z(x))$ defined in (\ref{z(x,+infini
determina+ reflection)}), (\ref{TiH(z(x))})  and
\begin{eqnarray} \label{z(x,t_1)reflextion}
z(x,t_1(E))&=&\frac{i}{c}\int_{t_1(E)}^x((E-V(t))^2-m^2c^4)^\frac12dt.
\end{eqnarray}
On $]-\infty,t_1(E)[$ the functions $((E-V(t))^2-m^2c^4)^\frac12$
and $\ \Ti H(z(x))$ are in $\R^+$. Then, we have
$$\W(\omin^-,\om_d^+)=2e^{i\pi/4}(1+h\phi(h))\exp\{\frac1h(z(t_1(E),-\infty)-z(t_1(E),+\infty))\},$$
where, $z(t_1,+\infty),\ z(t_1(E),-\infty)$ are defined respectively
in (\ref{z(x,+infini determina+ reflection)}) and (\ref{z(x,infty
determination+)}).

\medskip
{\it Computation of $\W(\om_d^+,\omout^-)$:} As in 
the previous paragraph  we extend $\om_d^+$ and the square roots
written there from $\Om^+(E)$ to
$\Om^-(E)\subset\C\setminus\{\im(z)>0,\ \re(z)=t_{1}\}$. We can also
find a path $\ti\gamma_2$ from $+\infty(1-i\ti\de_2)$ to
$-\infty(1+i\ti\de_2)$ (for $\ti\de_2>0$) transverse to the Stokes
lines along which we can extend $\om_d^+$. If we denote by
$\om_d^{+,2}$ the extension of $\om_d^+$ along $\ti\gamma_2$, we
have
$$\om_d^{+,2}=\exp\{\frac{1}{h}(+z(x,t_1(E))-z(t_1(E),+\infty))\}\left(\begin{array}{c}
e^{-i\pi/4}\Ti H(z(x))\\
e^{-i\pi/4}{\Ti H(z(x))^{-1}}
\end{array}\right)(1+h\phi(h)),$$
with $z(x,t_1(E))$, $z(t_1(E),+\infty)$ and $\Ti H(z(x))$ defined
respectively in $(\ref{z(x,t_1)reflextion})$, $(\ref{z(x,+infini
determina+ reflection)})$ and (\ref{TiH(z(x))}). On
$]-\infty,t_1(E)[$ the quantities $((E-V(t))^2-m^2c^4)^\frac12$ and
$ \Ti H(z(x))$ are in $\R^+$. Then, we have

$$\W(\om_d^+,\omout^-)=-2e^{-i\pi/4}(1+h\phi(h))\exp\{\frac1h(-z(t_1(E),-\infty)-z(t_1,+\infty))\},$$
where, $z(t_1,+\infty),\ z(t_1(E),-\infty)$ are defined respectively
in (\ref{z(x,+infini determina+ reflection)}) and (\ref{z(x,infty
determination+)}).

\medskip
{\it Computation of $\W(\omin^-,\omout^-)$:} Since the two solutions
$\omin^-,\omout^-$ are defined in $\Om^-(E)$, we compute the
Wronskian between these solutions as in (\ref{wronskianequal-2}) and
obtain
$$\W(\omin^-,\omout^-)=-2.$$

Then, we have
\begin{eqnarray*}
\al_{\rm out}^-&=&\frac{\W(\omin^-,\om_d^+)}{\W(\om_d^+,\omout^-)}
=-i(1+h\phi(h))\exp\{\frac2hz(t_1(E),-\infty)\}\\
\be_d^+&=&\frac{\W(\omin^-,\omout^-)}{\W(\om_d^+,\omout^-)}
=e^{i\pi/4}(1+h\phi(h))\exp\{\frac1h(z(t_1(E),-\infty)+z(t_1,+\infty))\},
\end{eqnarray*}
 with
 \begin{eqnarray*}
z(t_1(E),-\infty)&=& \frac{i}{c}\int_{-\infty}^{t_1(E)}\sqrt{(E-V(t))^2-m^2c^4}-\sqrt{(E-V^-)^2-m^2c^4}\ dt\label{z(t_1,-infini)t1}\\
\nonumber&+&\frac{i}{c}\sqrt{(E-V^-)^2-m^2c^4}\ t_1(E)\\
z(t_1(E),+\infty)&=& \frac1c\int_{+\infty}^{t_1(E)}\sqrt{m^2c^4-(E-V(t))^2}-\sqrt{m^2c^4-(E-V^+)^2}\ dt\label{z(t_1,+infini)t1}\\
\nonumber&+&\frac1c\sqrt{m^2c^4-(E-V^+)^2}\ t_1(E).
\end{eqnarray*}
This ends the proof of Theorem \ref{thmt1pt}.\\
\section{Zero mass case}\label{sec-mass0}

We suppose that $m=0$, $E\in]V^-,V^+[$ and $V$ satisfies assumption
(A),(see Fig. 11).
{\begin{center}
\psfrag{t1}[h][1][1][1]{$\!\!\!\!t_0(E)$}
\psfrag{Vl}[][][1][0]{$V^-$} \psfrag{Vr}[][][1][0]{$V^+$}
\psfrag{E}[][][1][0]{E} \psfrag{titre}[][][1][0]{Fig. 11. The double
turning point $t_0$}
\includegraphics[height=4cm, width=12cm]{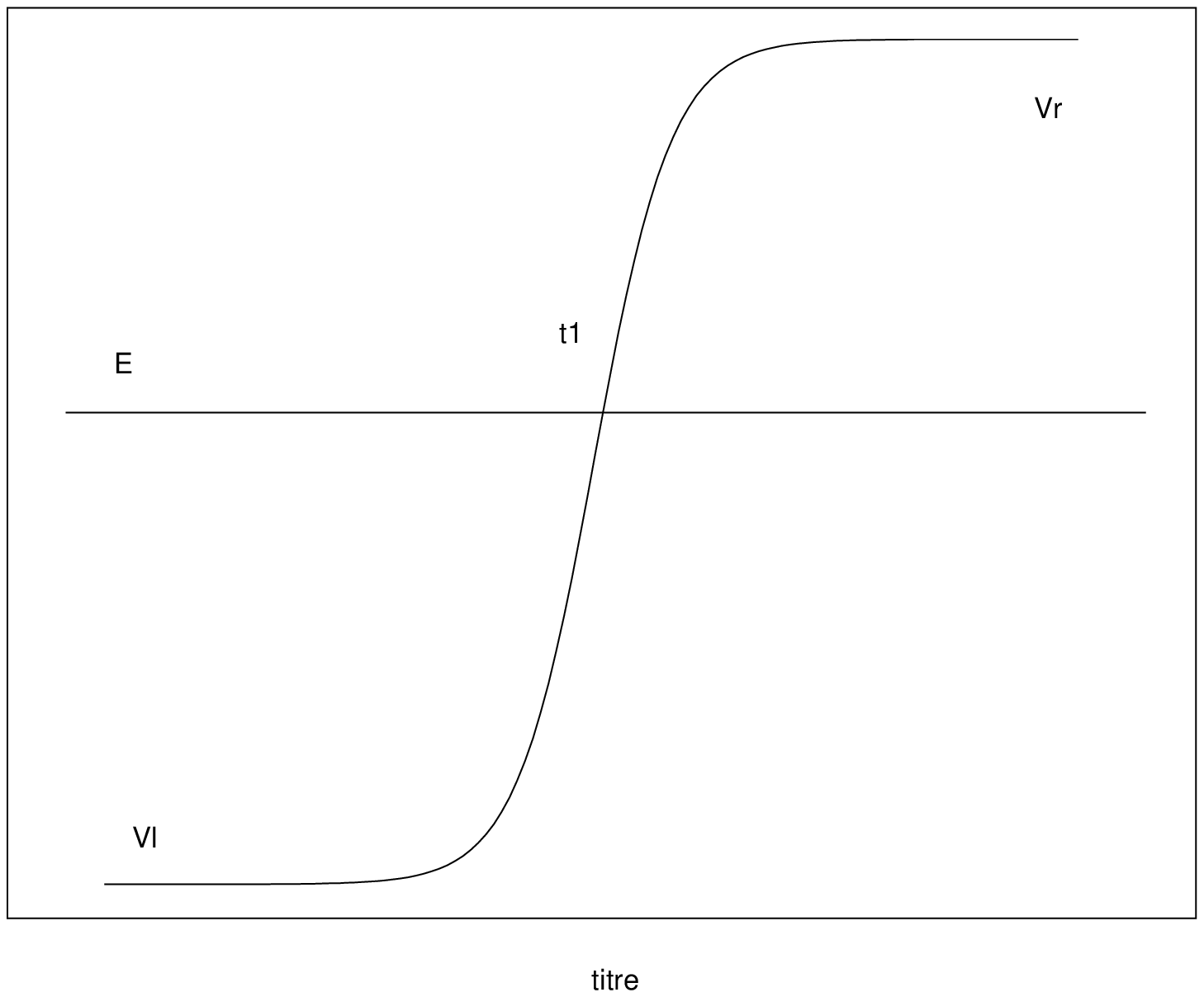}
\end{center}}
In this section we will work in two unbounded, simply-connected
domains $\Om_0^-(E)$, $\Om_0^+(E)$, where $\re V(x)<E$, $\re V(x)>E$
respectively and which coincide with $\calS$ for $|\re x|$
sufficiently large.
 Repeating the constructions of the solutions (see Sections \ref{sectionconst}, \ref{sectionJost}) and using Proposition \ref{propexpansion}
  for $m=0$,
there are two Jost solutions in
 $\Om_0^\pm(E)$:
 \begin{equation}\label{winout+-massnull}
\begin{aligned}
 \omin^{\pm}=&\exp\{\frac{1}{h}z_0(x,\pm\infty)\}\left(\begin{array}{c}
1\\
{\mp}1
\end{array}\right)(1+h\phi(h))
\\
 \sim&\exp\{\frac{\pm i}{hc}(V^\pm-E)x\}
\left(\begin{array}{c}
1\\
{\mp}1
\end{array}\right)\,\;\;\;\mbox{as }x\To\pm\infty,\\
 \omout^\pm=&\exp\{\frac{-1}{h}z_0(x,\pm\infty)\}\left(\begin{array}{c}
1\\
{\pm}1
\end{array}\right)(1+h\phi(h))\\
 \sim&\exp\{\frac{\mp i}{hc}(V^\pm-E)\
x\}\left(\begin{array}{c}
1\\
{\pm}1
\end{array}\right)\,\;\;\;\mbox{as }x\To\pm\infty,
\end{aligned}
 \end{equation}
with $\phi(h)$ a classical analytic symbol of non-negative order.
 The functions $z_0(x,\pm\infty)$ are defined by

\begin{eqnarray*}\label{eqchangeinfini0}
z_0(x,\pm\infty)&=&\pm\frac{i}{c}\int_{\pm\infty}^x(V(t)-V^\pm)dt\pm\frac{i}{c}(V^\pm-E)x.\\
\end{eqnarray*}

We suppose that there is only a simple zero of $V(x)-E$ (see Fig.  11). 
Recall that the turning points and the Stokes lines are also defined
for $m=0$ (see Definition \ref{defturningpoint}, Definition
\ref{defstokesline}). In our setting, there exists only a double
turning point $t_0(E)$ and the Stokes lines are described in Fig.
12.
{\begin{center}
\psfrag{x2-}[rh][][1][0]{}
\psfrag{t2}[][][1][0]{$t_2$}
\psfrag{x2+}[][][1][0]{}
\psfrag{t1}[][][1][0]{$\!\!\!t_0(E)$}
\psfrag{x1-}[][][1][0]{}
\psfrag{x1+}[][][1][0]{}
\psfrag{titre}[][][1][0]{Fig. 12. The Stokes lines and the double
turning point}
\includegraphics[height=5cm, width=14cm]
{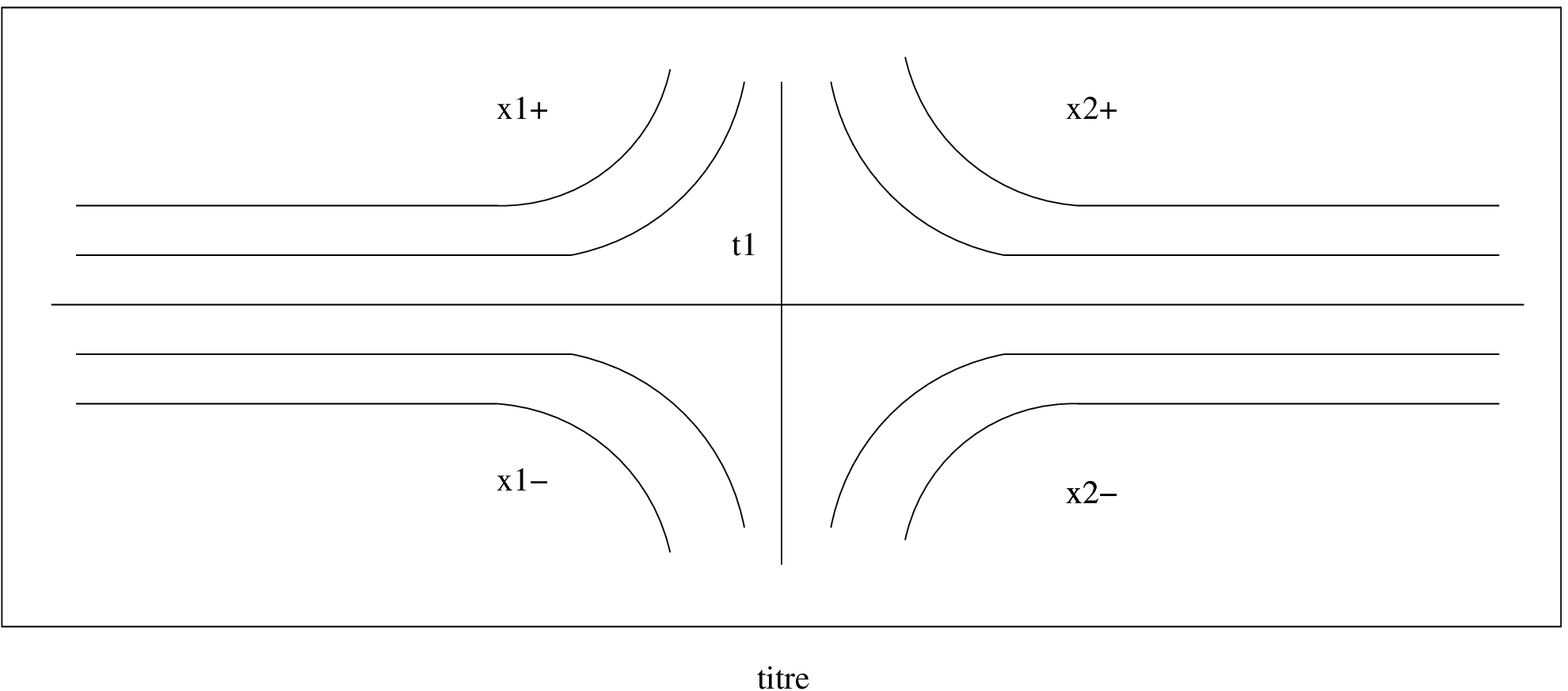}
\end{center}
}

In this section we use the definition of the scattering and transfer
matrix of Section \ref{sectionassuptionresutls} with the incoming
and outgoing Jost solutions defined above. As in Section
\ref{sectionJost}, it is sufficient to calculate the two terms
$r(E,h),\ t(E,h)$ (see (\ref{eq1tr}), (\ref{eq2tr})) to obtain the
matrix $\bbS$.

\medskip
{\it Computation of $\W(\omout^+,\omin^+)$:} Since the two solutions
$\omout^+,\,\omin^+$ are defined in $\Om_0^+$, we compute this
Wronskian in $\Om_0^+$ as (\ref{wronskianequal-2}) and we obtain
$$\W(\omout^+,\omin^+)=-2.$$

\medskip
{\it Computation of $\W(\omin^-,\omin^+)$:} The two solutions
$\omin^-,\omin^+$ are defined in $\Om_0^-(E),\Om_0^+(E)$
respectively. Since the Wronkian $\W(\omin^-,\omin^+)(x)$ is
independent on $x$ (see Remark \ref{rem wronsk ind x}), we compute
it in
$\Om_0^-(E)$ for example. For that, we extend $\omin^+$ from 
$\Om_0^+(E)$ 
to $\Om_0^-(E)$. Using (\ref{winout+-massnull}) we get
\begin{eqnarray*}
\W(\omin^-,\omin^+)&=& -2(1+h\phi(h)) \exp\left(\frac1h\left(z(x,-\infty)+z(x,+\infty)\right)\right),\\
&=& -2(1+h\phi(h))
\exp\left(\frac1h\left(z(t_0(E),-\infty)+z(t_0(E),+\infty)\right)\right).
\end{eqnarray*}

Then, we have
\begin{eqnarray}
t(E,h)&=&\frac{\W(\omin^-,\omin^+)}{\W(\omout^+,\omin^+)},\nonumber\\\label{eqvaleurt0(E,h)}
&=&(1+h\phi(h))\exp\left(\frac1h\left(z(t_0(E),-\infty)+z(t_0(E),+\infty)\right)\right). 
\end{eqnarray}
From (\ref{matrixS}), it follows
\begin{eqnarray*}
s_{11}&=&\frac{1}{\overline{\,t\,}(E,h)}=(1+h\phi_0(h))\exp\{\frac{i}{h}T_0(E)\},
\end{eqnarray*}
with
$$T_0(E)=\frac{1}{c}\left(\int_{+\infty}^{t_0(E)}(V(t)-V^+)dt-\int_{-\infty}^{t_0(E)}(V(t)-V^-)dt+t_0(E)(V^+-V^-)\right).
$$

From the form of the Stokes lines (see Fig. 12) there exist no path
transverse to the Stokes lines along which we calculate the Wronkian
$\W(\omout^-,\omin^+) $.
More precisely, the WKB method does not give the asymptotic behavior of
the Wronkian $\W(\omout^-,\omin^+) $ as well as the term
$$r(E,h)=\frac{\W(\omout^-,\omin^+)}{\W(\omout^+,\omin^+)}.$$

Nevertheless, (\ref{eqvaleurt0(E,h)}) together with the relation
(\ref{det=1}), implies:
\begin{eqnarray}\label{eqvaleurr}
r(E,h)=O(h).
\end{eqnarray}
Consequently, using (\ref{matrixS}), 
 we have
$$s_{21}=\frac{\overline{\,r\,}(E,h)}{\overline{\,t\,}(E,h)}=O(h)
,$$
$$
s_{12}=\frac{-r(E,h)}{\overline{\,t\,}(E,h)}=O(h)
.
$$
 This end the proof of Theorem \ref{thmmassnull}.

\appendix
\section{Spectrum of the Dirac operator
}\label{spectrumofdiracoperator}

\begin{prop}
Suppose that $V(x)$ is a $L^\infty$ application with values in the
space of Hermitian $2\time2-$matrix. Moreover we assume that
$$\|V(x)-V^\pm I_2\|\to0,\ \ \ \mbox{as}\ x\to\pm\infty.$$
Then the operator $H=H_0+V$ is a selfadjoint operator on $D(H_0)$
and
\begin{eqnarray}\label{sepectreess}
\si_{ess}(H)&=&
]-\infty,-mc^2+V^+]\cup[mc^2+V^-,+\infty[.\end{eqnarray}
\end{prop}
\begin{proof} In order to prove this proposition, we first calculate the essential spectrum of $H_0+W$, where $W$ is a $L^\infty$
potential with $W(x)=V^\pm I_2$ for $\pm x>R>0$
. From Lemma 5.1 of \cite{SRPB}, we know the following inclusion: 
\begin{eqnarray}
\si_{ess}(H_0+W)\subset
]-\infty,-mc^2+V^+]\cup[mc^2+V^-,+\infty[.\end{eqnarray}

Let us now prove the second inclusion. We denote
$$I^+:=]-\infty,-mc^2+V^+],\ \ \ \ \mbox{and }\ \ \ \ I^-:=[mc^2+V^-,+\infty[
.$$
For $E\in I^\pm$,  we consider the sequence
$$f_n^\pm=\exp\{\mp\frac{i}{hc}\Phi(E-V^\pm)
x\}\left(\begin{array}{c}
A(E-V^\pm)\\
\mp{A(E-V^\pm)^{-1}}
\end{array}\right)\chi(\pm x/n) \ \mbox{  for  }\ n\in\N,$$
with $\Phi(E)=\emph{sgn}(E)\sqrt{E^2-m^2c^4}$, $A(E)=\sqrt[4]{\frac{E+mc^2}{E-mc^2}}$ and $\emph{sgn}(E)=\frac{E}{|E|}$ for $E\not\in[-mc^2,mc^2]$. The function  $\chi\in C_0^\infty(\R)$ is such that $\chi(x)=1$ if $2<x<3$ and $\chi(x)=0$ if $x<1$ and $x>4$.

The normed sequence $(\frac{f_n^\pm}{\|f_n^\pm\|})_{n\in\N}$ has no
convergent subsequence and satisfies
$$(H_0+W-E)\frac{f_n^\pm}{\|f_n^\pm\|}\To0,\ \ \ \mbox{as}\ n\To+\infty.$$
 From the Weyl criterion, we deduce $I^\pm\subset\si_{ess}(H_0+W)$.
Consequently$$\si_{ess}(H_0+W)=
]-\infty,-mc^2+V^+]\cup[mc^2+V^-,+\infty[.$$ Finally, using Weyl's
theorem we obtain
$$\si_{ess}(H)=\si_{ess}(H_0+W),$$
and the proposition holds.
\end{proof}

\emph{Acknowledgments}. The author is grateful to V. Bruneau and
J.-F. Bony for many helpful discussions. We also thank the French
ANR (Grant no. JC0546063) for the financial support.

\end{document}